\documentclass{article}

\usepackage{arxiv}

\RequirePackage{amsthm,amsmath,amsfonts,amssymb}
\RequirePackage[numbers]{natbib}
\usepackage{amsmath}
\usepackage{color, graphicx}
\usepackage[colorinlistoftodos]{todonotes}
\usepackage{xcolor}
\usepackage{hyperref}
\usepackage{verbatim}
\usepackage{color,soul}
\usepackage{subfigure}
\usepackage{amsthm}
\usepackage{amsmath}
\usepackage{tikz}
\usepackage{float}
\usepackage{amsmath}
\usepackage{float}
\usepackage{siunitx}
\usepackage{amssymb}
\usepackage{comment}
\usepackage{multido}
\usepackage{graphicx}
\usepackage{svg}
\usepackage{marvosym}

\usepackage{wrapfig}
\usepackage{hyperref}
\usepackage[toc,page]{appendix}
\usepackage{bbm}
\usepackage{graphicx}
\usepackage{caption}

\usepackage{amsbsy}
\usepackage{lipsum}
\usepackage{multicol}
\usepackage{amsthm}
\usepackage{pifont}
\usepackage{bbm}
\usepackage{bm}
\usepackage{dsfont}
\usepackage{bigints}
\usepackage{mathtools}
\usepackage{graphicx}
\usepackage{svg}
\usepackage{calc}
\usepackage{scalerel}
\usepackage{accents}

\usepackage{mwe}
\usepackage{graphicx,wrapfig,lipsum}

\usepackage{mathtools}
\usepackage{mathrsfs}
\usetikzlibrary{mindmap,backgrounds}
\usepackage{blindtext}
\usepackage{hyperref}
\usepackage{pdfpages}
\usepackage[dvisvgm]{animate}
\usepackage[ruled,vlined]{algorithm2e}
\allowdisplaybreaks

\RequirePackage{graphicx}%% uncomment this for including figures

\usepackage[utf8]{inputenc}
\usepackage{amsmath,amsfonts,amssymb,amsthm}

\usepackage{enumitem}
\usepackage{letltxmacro}
\usepackage{nameref,hyperref}
\usepackage[capitalize]{cleveref}
\usepackage{indentfirst}

\newlist{thmlist}{enumerate}{1}
\setlist[thmlist]{label=\alph{thmlisti}., ref=\thetheorem.\alph{thmlisti},noitemsep}

\newlist{deflist}{enumerate}{1}
\setlist[deflist]{label=\alph{deflisti}., ref=\thedefinition.\alph{deflisti}.,noitemsep}

\numberwithin{equation}{section}
\newcommand{\pb}[0]{\mathbb{P}}
\newcommand{\ex}[0]{\mathbb{E}}
\newcommand{\ind}[0]{\mathbf{1}}
\newcommand{\ubar}[1]{\underaccent{\bar}{#1}}

\DeclareMathOperator*{\argmax}{arg\,max}
\DeclareMathOperator*{\argmin}{arg\,min}

\newlength\myindent
\newcommand\numberthis{\addtocounter{equation}{1}\tag{\theequation}}
\theoremstyle{plain}
\newtheorem{theorem}{Theorem}%[section]

\newtheorem{lemma}{Lemma}
\newtheorem{proposition}{Proposition}
\newtheorem*{definition*}{Definition}

\theoremstyle{remark}

\newtheorem{assumption}{Assumption}
\newtheorem{remark}{Remark}

\Crefname{theorem}{Theorem}{Theorems}

\title{Asymptotically efficient estimation under local constraint in Wicksell's problem}
\newcommand{\university}{\textit{Delft University of Technology, Mekelweg 4, Delft 2628CD, The Netherlands.}}
\newcommand{\authname}[1]{\textbf{#1}}
\author{
    \authname{Francesco Gili} \thanks{ \texttt{F.Gili@tudelft.nl}} \quad \quad \quad \authname{Geurt Jongbloed} \thanks{\texttt{G.Jongbloed@tudelft.nl}} \quad  \quad \quad \authname{Aad van der Vaart} \thanks{\texttt{A.W.vanderVaart@tudelft.nl}} \\
    \\
    \university % The university is placed below the authors
}
\hypersetup{
pdftitle={Asymptotically efficient estimation under local constraint in Wicksell's problem - F. Gili, G. Jongbloed, A. van der Vaart},
pdfsubject={},
pdfauthor={Francesco ~Gili, Geurt ~Jongbloed, Aad ~van der Vaart},
pdfkeywords={Wicksell's problem, Isotonic estimation, Efficiency theory},
}

\begin{document}
\maketitle

\begin{abstract}
We consider nonparametric estimation of the distribution function $F$ of squared sphere radii in the classical Wicksell problem. Under smoothness conditions on $F$ in a neighborhood of $x$, in \cite{21} it is shown that the Isotonic Inverse Estimator (IIE) is asymptotically efficient and attains rate of convergence $\sqrt{n / \log n}$. If $F$ is constant on an interval containing $x$, the optimal rate of convergence increases to $\sqrt{n}$ and the IIE attains this rate adaptively, i.e.\ without explicitly using the knowledge of local constancy. However, in this case, the asymptotic distribution is not normal. In this paper, we introduce three \textit{informed} projection-type estimators of $F$, which use knowledge on the interval of constancy and show these are all asymptotically equivalent and normal. Furthermore, we establish a local asymptotic minimax lower bound in this setting, proving that the three \textit{informed} estimators are asymptotically efficient and a convolution result showing that the IIE is not efficient. We also derive the asymptotic distribution of the difference of the IIE with the efficient estimators, demonstrating that the IIE is \textit{not} asymptotically equivalent to the \textit{informed} estimators. Through a simulation study, we provide evidence that the performance of the IIE closely resembles that of its competitors. 
\end{abstract}

% keywords can be removed

\section{Introduction}

In the field of stereology, scientists study the three-dimensional properties of materials and objects by interpreting their two-dimensional cross-sections. This may allow to estimate three-dimensional quantities without the use of expensive 3D reconstructions. In the Wicksell problem, a number of spheres are embedded in an opaque three-dimensional medium. Because the medium is opaque, we are not able to observe the spheres directly. However, we can observe a cross-section of the medium, which shows the circular sections of the spheres that happen to be cut by the plane. It is assumed that the spheres' squared radii are realizations from a cumulative distribution function (cdf) $F$, the object to be estimated. Following the same notation as in \cite{1,21}, a version of the density $g$ of the observable squared circle radii $Z$ is given by (cf.\ \cite{15}):
\begin{align}\label{eq: density of the observations}
g(z)=\frac{1}{2 m_0} \int_{z}^{\infty} \frac{d F(s)}{\sqrt{s-z}},    
\end{align}
where $ 0 < m_0=\int_0^{\infty} \sqrt{s} \, d F(s) < \infty$ is the expected sphere radius under $F$. Wicksell \cite{14} inverted this equation, by recognizing an Abel-type integral, and found an expression for $F$ in terms of $g$, given by:
\begin{align}\label{eq: inv}
    F(x)=1-\frac{\int_x^{\infty}(z-x)^{-1 / 2} g(z) \: dz}{\int_0^{\infty} z^{-1 / 2} g(z) \: dz} = 1 - \frac{V(x)}{V(0)}, \quad x>0,
\end{align}
where
\begin{align}\label{eq: V}
    V(x) \vcentcolon=\int_x^{\infty} \frac{g(z)}{\sqrt{z-x}} \: d z.
\end{align}
(Furthermore $F(x) = 0$ for all $x \leq 0$). Therefore, in order to estimate $F$ at a point $x>0$, the essential object to be estimated is the function $V$ at $x$ and at $0$. For a more detailed introduction to Wicksell's problem, we refer the reader to \cite{1,9,21}, where the asymptotic analysis of the reconstruction of $F$ from $V$ is also illustrated.

In this paper, we focus on the case in which the underlying cdf $F$ (and therefore $V$) is known to be constant on an interval of positive length. We conduct the analysis on a single interval $[\ubar{x},\bar{x}]$, for $\ubar{x}<\bar{x}$, of constancy. The same procedures can be applied if the underlying cdf $F$ is constant on multiple intervals. Note that even if $F$ is constant on  $[\ubar{x},\bar{x}]$, we do not observe any \textit{gap} in the data. Indeed, the density $g$ in \eqref{eq: density of the observations} is positive throughout the set $\{ x \, : \, F(x) < 1 \}$. In this setting, we compare an estimator (the "Isotonic Inverse Estimator" IIE, cf.\ \cite{1,21} for a detailed introduction) that attains the $\sqrt{n}$-rate of convergence adaptively, \textit{without} using the information that $V$ is constant on $[\ubar{x},\bar{x}]$, to three other \textit{informed} estimators that \textit{do use} this information and are therefore constrained to be constant on this interval. 

In Section \ref{sec: asymp distr}, we then show that all the informed estimators are asymptotically equivalent and normal, but not asymptotically equivalent to the IIE. Next in Section \ref{sec: local minimax}, we derive a lower bound for the local asymptotic minimax risk of any estimator sequence in this setting, and show that this is attained by the informed estimators. In contrast, the IIE is inefficient. In fact, we show that its (non Gaussian) limit distribution is the convolution of the limit distribution of the informed estimators and a nondegenerate factor. However, the difference between the estimators is small, as we illustrate in a simulation study, in Section \ref{sec:simulation study}.

A summary of the findings is that if it is known that the distribution function $F$ is constant on one or more intervals, then any of the informed estimators 
should be used for estimation on such intervals, resulting in efficient normal asymptotic behavior at $\sqrt n$-rate. Whenever this 
information is not available, it is preferable to use the isotonic estimator, which behaves closely enough to the informed estimators.

\subsubsection*{Motivation, connections with the literature and conclusions.} 

There is a vast literature on Wicksell's problem, see for instance \cite{1,2,3,9,10, 15, 16, 17, 18, 19, 20, 21}. Wicksell's problem has present-day applications as the estimation of the distribution of stars in a galactic cluster (cf.\ \cite{19}) or the estimation of the 3D microstructure of materials (cf.\ \cite{29,30}). In applications, it is common to see the use of the so-called Saltykov methods, based on numerical discretization. These methods are far from being efficient and in some cases not even consistent. The estimation of $F$ given data from $g$ is of interest not only for stereological procedures, but also from the mathematical point of view, both because of the unusual rates of convergence and because of the non-standard efficiency theory.

For nonparametric estimation the best attainable rate of convergence for estimation of the cdf $F$ at a point is $\sqrt{n/\log n}$. Many authors
contributed to this problem (\cite{1, 3, 10, 15, 16, 17, 18, 19, 20, 21}), introducing various estimators. We believe the state-of-the-art 
is the Isotonic Inverse Estimator, introduced in \cite{1}. (This regularizes the naive plug-in estimator obtained by replacing $g(z)\,dz$ in (1.3) by the
empirical distrubution of the inverse radii.) As proved in \cite{21}, this estimator maintains the rate of convergence $\sqrt{n/\log n}$ across
a range of smoothness conditions on $F$, with varying asymptotic variance,  and in this sense automatically adapts to smoothness, without
needing the selection of a bandwidth parameter. 

The case that $F$ is constant in an interval, treated in the present paper, can be viewed as the extreme case of smoothness in the sense of \cite{21}. The setting is of interest for practical 
applications in which the radii of the spheres are known to be restricted to certain values. The setting is also of theoretical interest, as the Isotonic Inverse Estimator is
not asymptotically normally distributed in this case.

\section{Construction of the estimators}\label{sec: construction estimators}

In this section, we introduce three projection-type \textit{informed} estimators \eqref{eq: proj isotonic smaller space}-\eqref{eq: proj naive smaller space}. We provide explicit constructions and prove that they solve determined minimization problems. 

Relation \eqref{eq: V} suggests a natural naive (empirical plug-in) estimator for $V$. For $\mathbb{G}_n$ the cdf of a sample $Z_1, \ldots, Z_n$ from the density $g$, let:
\begin{align}\label{eq: naive estimator of V}
    V_n(x) \vcentcolon=\int_x^{\infty} \frac{d\mathbb{G}_n(z)}{\sqrt{z-x}},
\end{align}
Using the "naive estimator", we construct all the other estimators as repeated projections of $V_n$ into designated spaces. Let $\ubar{x}<\bar{x}$, $V_n(x) = 0$ for $x>M$, and a fixed partition of $[0,M]$ into intervals. Let $\mathscr{P} = \{[x_i,x_{i+1})  \, : \, i=1,\ldots,I \}$ be such that $0 = x_0 < x_1 < x_2 < \cdots < x_I = M$, and $ [\ubar{x},\bar{x}] = [x_i,x_{i+1}]$ for some $i$. Define the spaces:
\begin{align*}
    &\mathcal{V} := \bigg\{ V : [0,\infty) \mapsto [0,\infty) \: : \: V \: \text{ decreasing, right-cont. }\bigg\}, \numberthis \label{eq: V space} \\
    &\overline{\mathcal{V}}_{\scriptstyle{\ubar{x},\bar{x}}} := \bigg\{ V : [0,\infty) \mapsto [0,\infty) \: : \: V \: \text{right-cont., piecewise-const. on } \mathscr{P} \bigg\}, \numberthis \label{eq: V space constant non monotonous} \\
    &\mathcal{V}_{\scriptstyle{\ubar{x},\bar{x}}} := \bigg\{ V \in \mathcal{V} \: : \: V \: \text{constant on } [\ubar{x},\bar{x}] \bigg\}. \numberthis \label{eq: V space constant} 
\end{align*}
All the above-defined spaces are convex cones. We define projections relative to the discrepancy measure $Q^{f} \: : \: \mathbb{L}_2 [0,\infty) \mapsto \mathbb{R}$, for a fixed function $f \in \mathbb{L}_1[0,\infty)$:
\begin{align}\label{eq: Q operator}
    h \, \mapsto \, Q^{f}(h) := \int_{0}^{\infty} h(s) \left(h(s) -2 f(s) \right) \, ds
\end{align}
Since $\|f-h\|_2^2=Q^f(h)+\|f\|_2^2$, minimising $h \rightarrow Q^f(h)$ is equivalent to minimising the $\mathbb{L}_2$-distance to $f$ in the case that $f$ is square-integrable. We use the discrepancy $Q^f$, because in our context the naive estimator $V_n$ is \textit{not} square-integrable. We consider as estimators the projections of the naive estimator $V_n$ onto the three spaces \eqref{eq: V space}-\eqref{eq: V space constant} and its repeated projection onto $\mathcal{V}$ and next on $\mathcal{V}_{\scriptstyle{\ubar{x},\bar{x}}}$:
\begin{align*}
  &\hat{V}_n = \underset{V \in \mathcal{V}}{\argmin} \: Q^{V_n}(V),  \numberthis \label{eq: isotonic estimator as proj} \\
   &\hat{V}^{\Pi}_n = \underset{V \in \mathcal{V}_{\scriptstyle{\ubar{x},\bar{x}}}}{\argmin} \: Q^{\hat{V}_n}(V).    \numberthis \label{eq: proj isotonic smaller space} \\
  & V^{\scriptscriptstyle{\mathscr{P}}}_n = \underset{V \in \overline{\mathcal{V}}_{\scriptstyle{\ubar{x},\bar{x}}}}{\argmin} \: Q^{V_n}(V),    \numberthis \label{eq: empirical slope as proj} \\
 & V^{\Pi}_n = \underset{V \in \mathcal{V}_{\scriptstyle{\ubar{x},\bar{x}}}}{\argmin} \: Q^{V_n}(V), \numberthis \label{eq: proj naive smaller space} 
\end{align*}
(We note that projecting on a cone and then projecting again onto a subcone does not need to give the same result as projecting directly onto the subcone). Existence and uniqueness of the solutions of these minimization problems can be shown along the lines of Theorem 1.2.1 in \cite{28}. The first estimator $\hat{V}_n$ is the Isotonic Inverse Estimator (IIE) introduced in \cite{1}. The other three estimators are informed estimators that take the local constancy of $F$ into account.

The following proposition gives explicit constructions of the estimators \eqref{eq: isotonic estimator as proj}-\eqref{eq: proj naive smaller space} (the first statement of the following proposition has been proved in \cite{1}). For a function $k$ on $[0,M]$, the Least Concave Majorant (LCM) and its right-hand side derivative, are given by:
\begin{align*}
    &k^* := \argmin {\{ f \, : \, f(x) \geq k(x) \, \, \text{for} \, \, x \in [0,M], \, f \, \text{concave} \}}, \\
    & (k^*)^{\prime} (x):= \inf_{u <x} \sup_{v \geq x} \frac{k(v)-k(u)}{v-u}.
\end{align*}
By definition, $k^*(x) \geq k(x)$ for all $ x \in [0,M]$, and $k^*(0)=k(0)$, $k^*(M) = k(M)$. Define functions $U_n$ and $U$ as
\begin{align*}
    U_n(x)&=\int_0^x V_n(y) \, dy=2 \int_0^{\infty} \sqrt{z} \, d\mathbb{G}_n(z)-2 \int_x^{\infty} \sqrt{z-x} \, d \mathbb{G}_n(z), \numberthis \label{eq: U_n} \\
    U(x)&=\int_0^x V(y) \, d y =\frac{\pi}{2 m_0} \int_0^x(1-F(y)) \, d y. \numberthis \label{eq: def U}
\end{align*}
The last equality follows from \eqref{eq: inv}. 

\begin{proposition}[Construction estimators]\label{prop: construction estimators}
Let $U_n$ be as in \eqref{eq: U_n}, $\ubar{x}<\bar{x}$ with $[\ubar{x},\bar{x}) \in \mathscr{P}$ and $x \in [0,\infty)$. The solutions of the minimizations \eqref{eq: isotonic estimator as proj}-\eqref{eq: proj naive smaller space} are:
\begin{align*}
     &\hat{V}_n (x) = (U^*_n)^{\prime} (x). \numberthis \label{def: isotonic estimator} \\
     & \hat{V}^{\Pi}_n(x) = \begin{cases} 
            \frac{U^*_n(\bar{x}) - U^*_n(\ubar{x})}{\bar{x}-\ubar{x}}, \quad x \in [\ubar{x},\bar{x}], \\
            \hat{V}_n(x), \quad x \notin [\ubar{x},\bar{x}].
            \end{cases} \numberthis \label{eq: projection iso} \\
     & V^{\scriptscriptstyle{\mathscr{P}}}_n (x) = V^{(x_i,x_{i+1})}_n (x), \quad  x \in [x_i, x_{i+1}), \quad i = 1, 2, \ldots, I, \numberthis 
\end{align*}
where for $x \in [x_i,x_{i+1}) \in \mathscr{P}$:
\begin{align*}
    V^{(x_i,x_{i+1})}_n (x) =
            \frac{U_n(x_{i+1}) - U_n(x_i)}{x_{i+1}-x_{i}}, \numberthis \label{eq: natural competitor}
\end{align*}
\vspace{-0.3cm}
\begin{align*}
    & V^{\Pi}_n(x) = (\widetilde{U}^*_n)^{\prime} (x), \quad\quad\quad\quad\quad\quad\quad\quad\quad\quad \quad\quad\quad\quad\quad  \quad\, \,     \numberthis \label{eq: projection naive} 
\end{align*}
where:
\begin{align}\label{eq: U_n tilde}
\widetilde{U}_n(t) := \begin{cases}
    \frac{U_n(\bar{x})-U_n(\ubar{x})}{\bar{x}-\ubar{x}} (t-\ubar{x}) + U_n(\ubar{x}), \quad t \in [\ubar{x},\bar{x}], \\
    U_n(t), \quad t \notin [\ubar{x},\bar{x}].
\end{cases}
\end{align}
\end{proposition}

The \hyperlink{proof of prop constructions}{proof} of Proposition \ref{prop: construction estimators} is given in Appendix A and the asymptotic properties of the constructed estimators are given in Theorem \ref{thm: distrib flat F at x total} below. 

\begin{remark}
For any $v \in [\ubar{x},\bar{x}]$ and $x \in [\ubar{x},\bar{x}]$ with $v \neq x$, we have:
$$
U(v) - U(x) - V(x)(v-x) = 0,
$$
and thus by evaluating the above relation at $v=\bar{x}$ and $x=\ubar{x}$ we get:
\begin{align}\label{eq: V if F constant}
U(\bar{x}) - U(\ubar{x}) - V(x)(\bar{x}-\ubar{x}) = 0.
\end{align}
This gives an intuition for \eqref{eq: natural competitor} and explains why $V^{(\ubar{x},\bar{x})}_n(x)$, i.e. $V^{\mathscr{P}}_n (x)$ for $x \in [\ubar{x},\bar{x})$, can be called "empirical slope" on $[\ubar{x},\bar{x})$. It what follows we use the notation $V^{(\ubar{x},\bar{x})}_n(x)$ to indicate $V^{\mathscr{P}}_n (x)$ for $x \in [\ubar{x},\bar{x})$.
\end{remark}

The isotonic estimator in \eqref{eq: isotonic estimator as proj} can be implemented using \eqref{def: isotonic estimator}, where the least concave majorant of a function can be computed by classical algorithms like the PAVA algorithm (cf.\ \cite{2}). Similarly, the solutions to \eqref{eq: proj isotonic smaller space}-\eqref{eq: proj naive smaller space} can be easily implemented using \eqref{eq: projection iso}, \eqref{eq: natural competitor} and \eqref{eq: projection naive}. Alternatively, the solution to the minimization problem \eqref{eq: proj naive smaller space} can be computed algorithmically using a \textit{profile procedure} as follows. First, fix $a>0$ and compute the minimizer of $Q^{V_n}(V)$ over the space of decreasing, right-continuous functions $V$ that have value $a$ on $[\ubar{x},\bar{x}]$, and call such minimizer $V^a_n$. Next we optimize the mapping $a \mapsto Q^{V_n}(V^a_n)$ over $a>0$ and obtain the desired projection. We give a detailed version of \hyperlink{algo L2 proj}{this algorithm} at the end of the appendix.

\section{Asymptotic distributions of the estimators}\label{sec: asymp distr}

In \cite{21} (cf.\ Theorem 2), it was shown that if $V \in \mathcal{V}_{\ubar{x},\bar{x}}$ and $x \in (\ubar{x},\bar{x})$, then as $n \rightarrow \infty$:
    \begin{align}\label{eq: adaptiveness iso}
        \sqrt{n} \left( \hat{V}_n(x) - V(x) \right) \rightsquigarrow L_x,
    \end{align}
    where, for any $a \in \mathbb{R}$: $\pb ( L_x \leq a ) = \pb ( \argmax_{s \in [\ubar{x},\bar{x}]} \left\{ \mathbb{Z}_x(s) - a s\right\} \leq x ),$
    for $(x,t) \mapsto \mathbb{Z}_x(t)$ a centered continuous Gaussian Process with covariance structure given by:
    \begin{align}\label{eq: cov GP paper 1}
        \mathrm{Cov} \left( \mathbb{Z}_y(t),\mathbb{Z}_x(s) \right) = 4 \, \mathrm{Cov} \left(  \sqrt{  (\smash[b]{Z_y})_{+}} - \sqrt{(Z_t)_{+}} , \sqrt{(Z_x)_{+}} - \sqrt{(Z_s)_{+}}  \right).
    \end{align}
    where $(\cdot)_+ = \max{\{0, \cdot \}}$ and $Z_{u} := Z-u$, $Z \sim g$, $u \in \mathbb{R}$. In the present section, we use the characterizations of the \textit{informed} estimators to obtain the first main result of this paper that puts the different estimators in perspective. Assume that:
\begin{align*}
    \quad \int_{0}^{\infty} s^{\frac{3}{2}} \: dF(s)  < \infty. \numberthis \label{eq: finite first moment}
\end{align*}

\begin{theorem}[Asymptotics]\label{thm: distrib flat F at x total}
     Let $F$ be the distribution function of the squared sphere radii and $g$ the corresponding density of the squared circle radii $Z$ according to \eqref{eq: density of the observations}. Let $x \geq 0$, and $K:=[\ubar{x},\bar{x}]$, for $\ubar{x}<\bar{x}$, the biggest interval that contains $x$ on which $F$ is constant and let \eqref{eq: finite first moment} hold true. Then, as $n \rightarrow \infty$:
    \begin{thmlist}
        \item The estimators $V^{(\ubar{x},\bar{x})}_n$, $V^{\Pi}_n$ and $\hat{V}^{\Pi}_n$ are all asymptotically equivalent: \label{thm: distrib flat F at x total 2}
                \begin{align*}
                    \sqrt{n} \left(  V^{\Pi}_n(x) -V^{(\ubar{x},\bar{x})}_n (x) \right) = o_p(1), \quad \quad \sqrt{n} \left(  \hat{V}^{\Pi}_n(x) -V^{(\ubar{x},\bar{x})}_n (x) \right) = o_p(1).
                \end{align*}
                Furthermore, they all attain the same limiting distribution, given by:
                \begin{align*}
                    \sqrt{n} \left(  V^{(\ubar{x},\bar{x})}_n (x) - V (x) \right) \rightsquigarrow N(0,\sigma^2_{\ubar{x},\bar{x}})
                \end{align*}
                where: 
                \begin{align}\label{eq: sigma}
                    \sigma^2_{\ubar{x},\bar{x}} = \mathrm{Var}  \left( 2 \, (\bar{x}-\ubar{x})^{-1} \big\{ \sqrt{\smash[b]{(Z_{\ubar{x}})_{+}}} - \sqrt{\smash[b]{(Z_{\bar{x}})_{+}}} \big\} \right).
                \end{align}
                \vspace{-0.3cm}
        \item The isotonic estimator $\hat{V}_n$ is \textbf{not} asymptotically equivalent to $V^{(\ubar{x},\bar{x})}_n$, as: \label{thm: distrib flat F at x total 1}
            \begin{align}\label{eq: limiting W}
                \sqrt{n} \left( \hat{V}_n(x) -V^{(\ubar{x},\bar{x})}_n (x) \right) \rightsquigarrow W,
            \end{align}
            where, for any $a \in \mathbb{R}$, $\pb \left( W \leq a \right) = \pb \left( \argmax_{s \in K} \left\{ \mathbb{Z}(s) - a s\right\} \leq x \right) $, for
            \begin{align*}
           & \mathbb{Z}(t) = \frac{\bar{x}-t}{\bar{x}-\ubar{x}} \mathbb{Z}_{\ubar{x}}(t) + \frac{t-\ubar{x}}{\bar{x}-\ubar{x}} \mathbb{Z}_{\bar{x}}(t), \numberthis \label{eq: gauss process}
            \end{align*}
            and $(t,x) \mapsto \mathbb{Z}_x(t)$ the centered Gaussian Process with covariance structure \eqref{eq: cov GP paper 1}.
        \item Under the assumptions of Proposition \ref{prop: LAN}, the sequences $$\sqrt{n} \left(  V^{(\ubar{x},\bar{x})}_n (x) - V (x) \right) \, \text{and} \, \, \, \sqrt{n} \left( \hat{V}_n(x) -V^{(\ubar{x},\bar{x})}_n (x) \right), $$ are asymptotically independent. \label{thm: asymptotic indep}
    \end{thmlist}
\end{theorem}

The \hyperlink{proof of thm 1 a}{proofs} of Theorem \ref{thm: distrib flat F at x total 2}, \ref{thm: distrib flat F at x total 1} and \ref{thm: asymptotic indep} are given in Appendix A. Combining \ref{thm: distrib flat F at x total 2}, \ref{thm: distrib flat F at x total 1} and \ref{thm: asymptotic indep} of the Theorem, we see that the limit $L_x$ of the isotonic estimator in \eqref{eq: adaptiveness iso} is the convolution of the normal limit $N (0, \sigma^2_{\ubar{x},\bar{x}})$ of the informed estimators and the distribution of $W$ in \eqref{eq: limiting W}. Thus the limit distribution of the isotonic estimator is less concentrated than the limit distribution of the informed estimators.

\section{Lower bound for the local asymptotic minimax risk}\label{sec: local minimax}

In this section, we show that the estimator \eqref{eq: empirical slope as proj} is asymptotically efficient, in the sense that its asymptotic variance is the smallest attainable. This proves that also \eqref{eq: proj naive smaller space} and \eqref{eq: proj isotonic smaller space}, which are asymptotically equivalent to \eqref{eq: empirical slope as proj}, are efficient. The claim is valid in the sense of both the local asymptotic minimax theorem and the convolution theorem (c.f.\ Chapter 25 of \cite{6}). The key is a LAN expansion for a submodel constructed by perturbing the true function $V$ in the least favorable direction.

We consider the estimation of $V$ at a point $x$ in $(\ubar{x},\bar{x})$, knowing that $F$, and thus $V$, are constant on this interval (cf.\ \eqref{eq: def V space tangent}). 
The associated version of the density of the observations in Wicksell's problem is expressed in $V$ as:
\begin{align}\label{eq: density observations dep V}
g_{\scriptscriptstyle{V}}(z)  = - \int_{z}^{M} \frac{dV(s)}{\pi \sqrt{s-z}} \quad  \left( = \frac{1}{2m_0} \int_{z}^{M} \frac{dF(s)}{\sqrt{s - z}}  \right).
\end{align}
We denote the corresponding probability measure by $G_{\scriptscriptstyle{V}}$. We dexcribe a  so called "least favorable" submodel in terms of $V$. We assume that the true parameter $V$ satisfies the following assumption.
\begin{assumption}\label{assumption: LAM}
    Let $V$ be supported on $[0,M]$, constant on $[\ubar{x},\bar{x}]$, for $\ubar{x}<\bar{x}$ and possess Lipschitz continuous density $v$ on $[\ubar{x},\bar{x}]^c$. Moreover $\exists \, \, \eta >0$ such that:
    \begin{align}\label{eq: L^2 condition}
        \int_{[\ubar{x},\bar{x}]^c} \frac{1}{|v(s)|} \, ds +  \int_{\ubar{x}-\eta}^{\ubar{x}} \frac{\log^2{(\ubar{x}-s)}}{|v(s)|} \, ds + \int_{\bar{x}}^{\bar{x} + \eta} \frac{\log^2{(s-\bar{x})}}{|v(s)|} \, ds < \infty.
    \end{align}
\end{assumption}
Define functions $k^{(\ubar{x},\bar{x})}$, $h$ and $h_0$ by, for $\ubar{x} < \bar{x}$:
\begin{align*}
        &k^{(\ubar{x},\bar{x})} (z):= 2 \frac{\sqrt{(z-\ubar{x})_{+}} - \sqrt{(z-\bar{x})_{+}}}{\bar{x}-\ubar{x}},\\
        &h (z) := \begin{cases}
            - \frac{1}{\sqrt{\pi} v(z)} \bigintssss_{z}^{M} \frac{(k^{(\ubar{x},\bar{x})} g_{\scalebox{0.4}{$V$}} )^{\prime} (s)}{\sqrt{s-z}} \, ds,  \quad \quad &z \notin [\ubar{x},\bar{x}], \\
             0, \quad \quad &z \in [\ubar{x},\bar{x}],
        \end{cases} \numberthis \label{eq: h funtion}\\
        &h_0 (z) := h(z) + \frac{2}{\pi} \int_{0}^M \sqrt{s} h(s) \, dV(s). \numberthis \label{eq: h_0} 
\end{align*}
The \textit{least favorable submodel} is defined by:
\begin{align}\label{eq: actual perturbation}
    dV_t(s) = \frac{1}{c_t}(1+ t h_0(s))_{+} \, dV(s), \quad c_t := - \frac{2}{\pi} \int_{0}^{M} \sqrt{s} (1 + th_0(s))_{+} \, dV(s).
\end{align} 
Note that $ \int_0^{\infty} \sqrt{s} \, d V_t(s)=-\frac{\pi}{2}$, that $c_t \geq 0$ and that the obtained $V_t$ is decreasing. In Appendix B we show that the score function of the model $t \mapsto G_{\scriptscriptstyle{V_t}}$ at $t=0$ is the centered function $k^{(\ubar{x},\bar{x})}$, the influence function of the informed estimators. Assumption \ref{assumption: LAM} ensures that indeed $h_0$ is mapped to $k^{(\ubar{x},\bar{x})}$ via the score operator and that $- \frac{2}{\pi}\int h_0^2 (s) \sqrt{s} \, dV(s) < \infty$. A detailed explanation is contained in the study of the tangent spaces in the context of our problem in Appendix B. 

The next result shows the model $t \mapsto G_{\scriptscriptstyle{V_t}}$ is locally asymptotically normal (LAN, see \cite{6}, Chapter 7). This result is fundamental because both Theorem \ref{thm: asymptotic indep} and Theorem \ref{thm: LAM} rely on it. The \hyperlink{Proof of LAN prop}{proof} of Proposition \ref{prop: LAN} is given in Appendix A. 

\begin{proposition}[LAN expansion]\label{prop: LAN}
Let $V_{t}$ be as in \eqref{eq: actual perturbation} and $t_n = t/ \sqrt{n}$ with $t \in \mathbb{R}$. If \eqref{eq: finite first moment} and Assumption \ref{assumption: LAM} hold true, then:
\begin{align}\label{eq: linear expansion}
    \sum_{i=1}^n \log{\frac{g_{\scriptscriptstyle{V_{t_n}}}(Z_i)}{g_{\scriptscriptstyle{V}}(Z_i)}} = t \, \Delta_n - \frac{t^2}{2} \, \sigma^2_{\ubar{x},\bar{x}} + o_p(1)
\end{align}
where: $\Delta_n := \frac{1}{\sqrt{n}} \sum_{i=1}^n  \left( k^{(\ubar{x},\bar{x})}(Z_i) - \mathbb{E}_{\scriptscriptstyle{G_{\scriptscriptstyle{V}}}} k^{(\ubar{x},\bar{x})} \right) \rightsquigarrow N (0, \sigma^2_{\ubar{x},\bar{x}})$ and $\sigma^2_{\ubar{x},\bar{x}}$ is as in \eqref{eq: sigma}.
\end{proposition}

The second fundamental result of this paper gives a lower bound for the local asymptotic minimimax (LAM) risk of an arbitrary estimator, and the related convolution theorem. 

\begin{theorem}[LAM \& Convolution]\label{thm: LAM}
    Let \eqref{eq: finite first moment} and Assumption \ref{assumption: LAM} hold true and let $\sigma^2_{\ubar{x},\bar{x}}$ be as in \eqref{eq: sigma}. Let then $\ell \, : \, \mathbb{R} \mapsto [0,\infty)$ be symmetric and subconvex. Then for any $x \in (\ubar{x},\bar{x})$, and every estimator sequence $\left( V_n(x)\right)_{n \in \mathbb{N}}$:
    \begin{align}\label{eq: LAM}
    \sup_{\substack{I \subset \mathbb{R} \\ I \: \textup{finite}}} \liminf_{n \rightarrow \infty} \, \sup_{t \in I} \, \ex_{\scriptscriptstyle{G}_{\scriptscriptstyle{{t/\sqrt{n}}}}} \, \ell\left(\sqrt{n}\left(V_n(x)- V_{t/\sqrt{n}}(x)\right)\right) \geq \int \ell \, d N (0, \sigma^2_{\ubar{x},\bar{x}}),
    \end{align}
    Moreover, for any estimator sequence such that $\sqrt{n}\left(V_n(x)- V_{t/\sqrt{n}}(x)\right) \rightsquigarrow K$ under $G_{\scriptscriptstyle{t/\sqrt{n}}}$, for every $t \in \mathbb{R}$, and some given distribution $K$, this distribution is the convolution of $N (0, \sigma^2_{\ubar{x},\bar{x}})$ and some other distribution.
\end{theorem}
 
The \hyperlink{Proof of LAM thm}{proof} of Theorem \ref{thm: LAM} is given in Appendix A. By Theorem \ref{thm: distrib flat F at x total} the informed estimators attain equality in the LAM risk (for bounded continuous loss functions), and their limit distribution is $N (0, \sigma^2_{\ubar{x},\bar{x}})$ without extra convolution factor. On the other hand, the limit distribution of the isotonic estimator contains the distribution of $W$ as in \eqref{eq: limiting W} as extra factor. Hence the isotonic estimator is not LAM.

\section{Simulation study}\label{sec:simulation study}

In this section, we present a simulation study that shows that the isotonic inverse estimator, even if not efficient, behaves very closely to the efficient estimator. In practice, the information whether or not the cdf $F$ is constant on a specific interval is generally unavailable. In that situation the best choice would be to use the isotonic estimator, given the efficiency theory and the adaptivity results developed in \cite{21}. However, if the cdf is known to be constant on some known interval, then the isotonic estimator will incur a loss in terms of efficiency. In this section, we perform a simulation study to illustrate the difference in terms of performance between the isotonic estimator $\hat{V}_n$ and the efficient informed estimator $V^{(\ubar{x},\bar{x})}$. We consider an underlying distribution $F$ of the squared sphere radii, which is constant on $[2,3]$, namely:

\begin{figure}[H]
    \centering
    \vspace{-0.3cm} 
    \begin{minipage}[b]{0.45\textwidth} % Adjust the width as needed
        \centering
        \includegraphics[width=6.4cm]{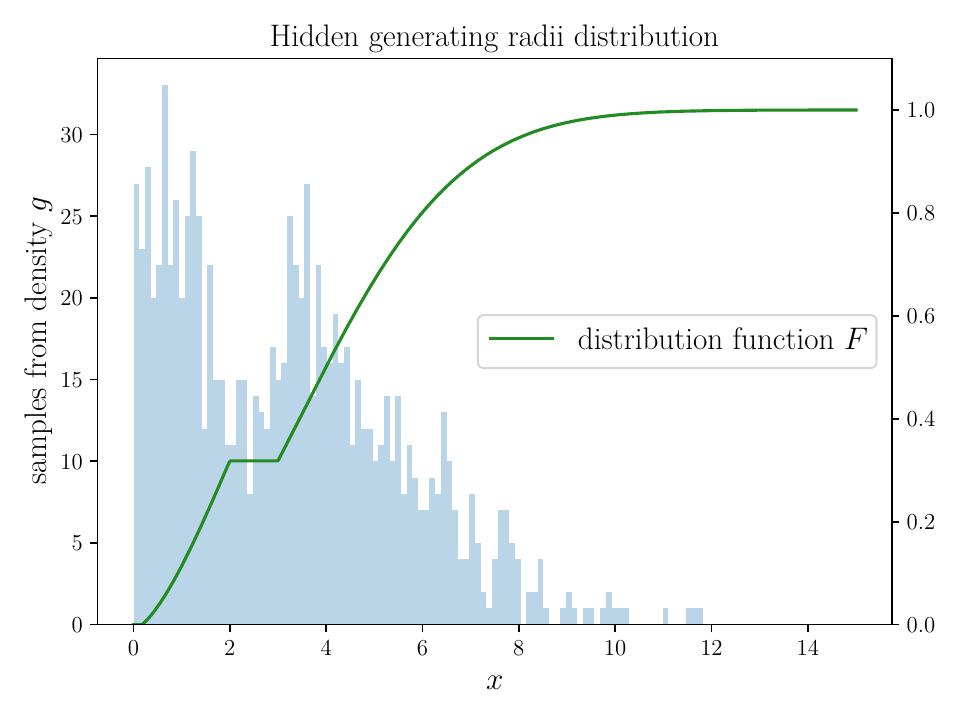}
        \label{fig:image11}
    \end{minipage}
    \hfill
    \begin{minipage}[b]{0.45\textwidth} % Adjust the width as needed
        \centering
        \raisebox{3.1cm}{ % Adjust the value here to move the formula up
        \resizebox{\columnwidth}{!}{%
            \begin{minipage}{\textwidth}
                \begin{align*}
                    &F(x) = \begin{cases}
                        1 - e^{-x^2/20}, & 0 \leq x < 2 \\
                        1 - e^{-1/5}, & x \in [2,3] \\
                        1 - e^{-(x-1)^2/20}, & x > 3
                    \end{cases}
                    \\
                    &g_{\scriptscriptstyle{V}}(x) = \frac{1}{2m_0} \int_{x}^{M} \frac{s e^{-\frac{s^2}{10}} \ind_{[0,2]}(s)}{10\sqrt{s - x}} \, ds \\
                    &\hspace{-0.2cm}+ \frac{1}{2m_0} \int_{x}^{M} \frac{(s-1)e^{-\frac{(s-1)^2}{20}} \ind_{[3,M]}(s)} {10\sqrt{s - x}} \, ds
                \end{align*}
            \end{minipage}}
        }
    \end{minipage}
    \vspace{-0.8cm} % Adjust the value here to decrease the space
    \caption{Underlying cdf $F$ and in light blue a histogram of a sample of size $1000$ from $g$.}\label{img: cdf1}
    \vspace{-0.5cm}
    \label{fig:both_images1}
\end{figure}

\begin{figure}[htbp]
    \centering
    \vspace{-0.1cm} 
    \begin{minipage}[b]{0.45\textwidth} % Adjust the width as needed
        \centering
        \includegraphics[width = 6.8cm]{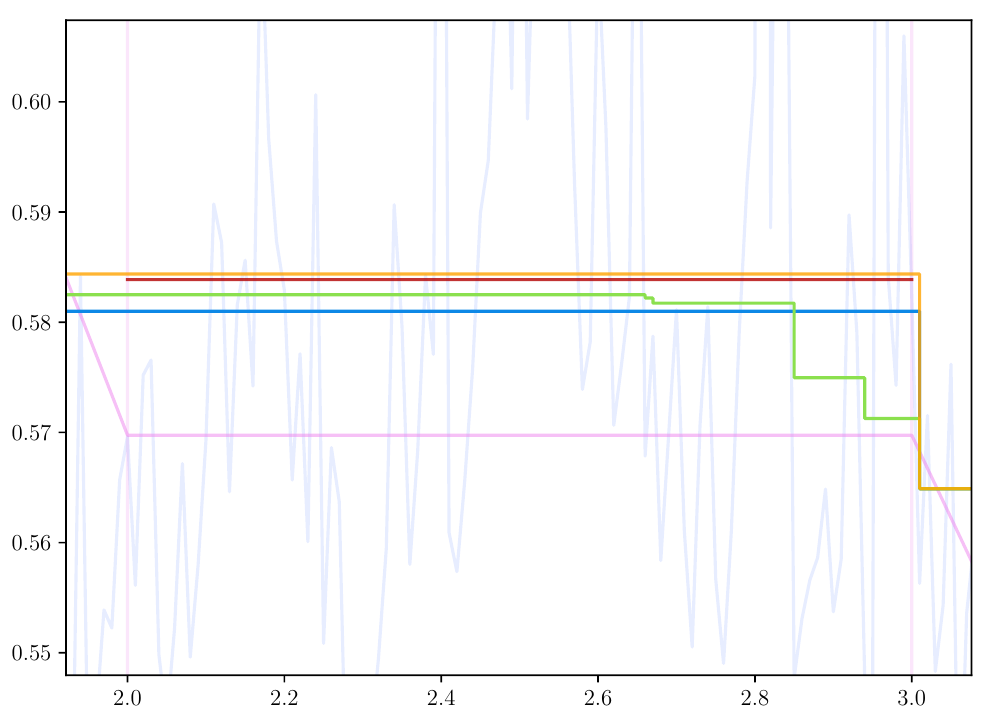}
        \label{fig:image1}
    \end{minipage}
    \hfill
    \begin{minipage}[b]{0.45\textwidth} % Adjust the width as needed
        \centering
        \includegraphics[width = 4.3cm]{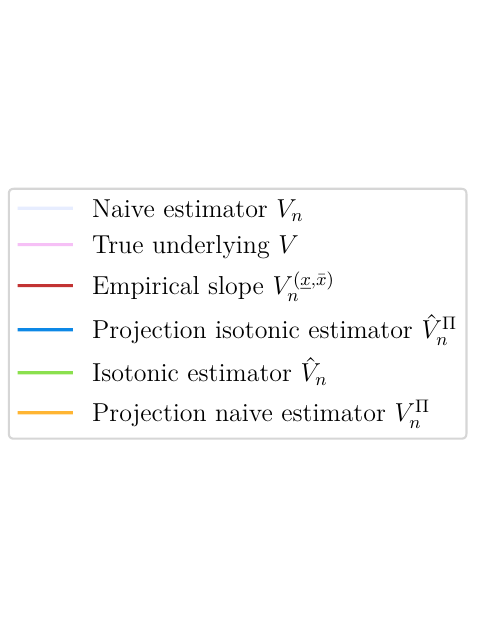}
        \label{fig:image2}
    \end{minipage}
    \vspace{-0.6cm} % Adjust the value here to decrease the space
    \caption{All estimators and the true function $V$, which constant on $[2,3]$, for a sample of size $n=300$.}\label{img: difference all estimators}
    \vspace{-0.1cm}
    \label{fig:both_images}
\end{figure}

Figure \ref{img: difference all estimators} shows the true underlying $V$ with all estimators defined in this paper within the specified interval. As $\hat{V}^{\Pi}_n$, $V^{\Pi}_n$, and $V^{(\ubar{x},\bar{x})}_n$ are all asymptotically equivalent by  Theorem \ref{thm: distrib flat F at x total}, we focus on the comparison between $V^{(\ubar{x},\bar{x})}_n$ and $\hat{V}_n$.

\begin{figure}[H]
    \centering
    \vspace{-0.7cm}
    \includegraphics[width=6.5cm]{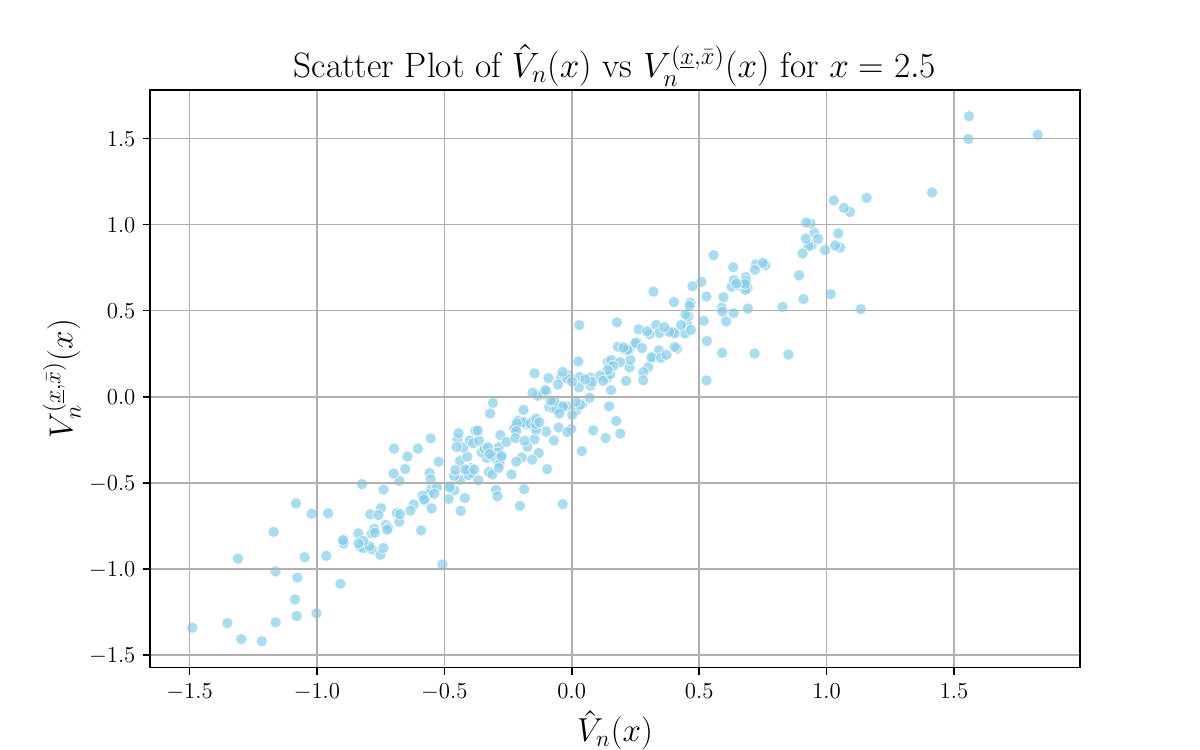}
    \hspace{-0.8cm}
    \includegraphics[width=6cm]{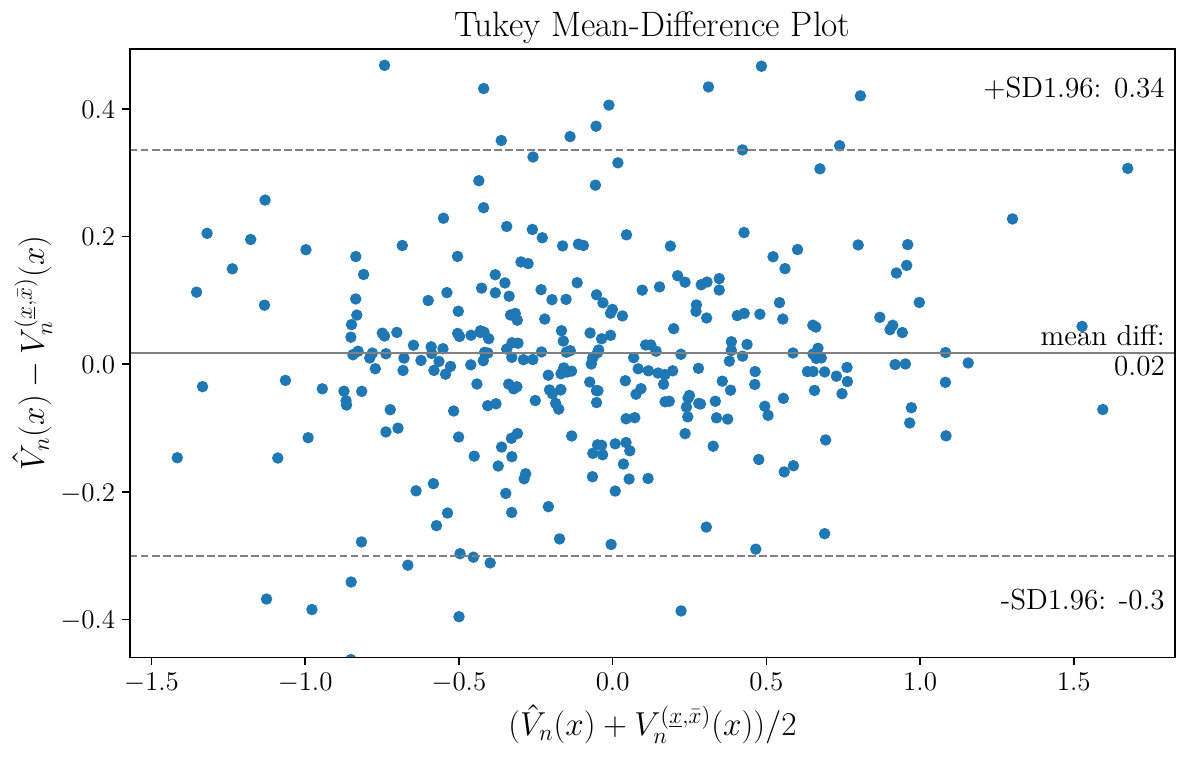}
    \setlength{\belowcaptionskip}{-6pt}
    \caption{Scatter and Tukey mean-difference plots for $n=1000$ and $300$ repetitions of $\hat{V}_n$ and $V^{(\underline{x},\overline{x})}_n.$}\label{img: scatter + tuckey}
  \end{figure}
\begin{figure}[H]
    \centering
    \vspace{-0.2cm}
    \includegraphics[width=12.1cm]{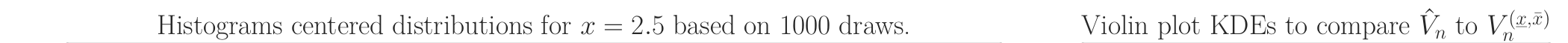}
    \vspace{-0.05cm} 
    
    \begin{minipage}[b]{\linewidth}
      \centering
      \includegraphics[width=12.1cm]{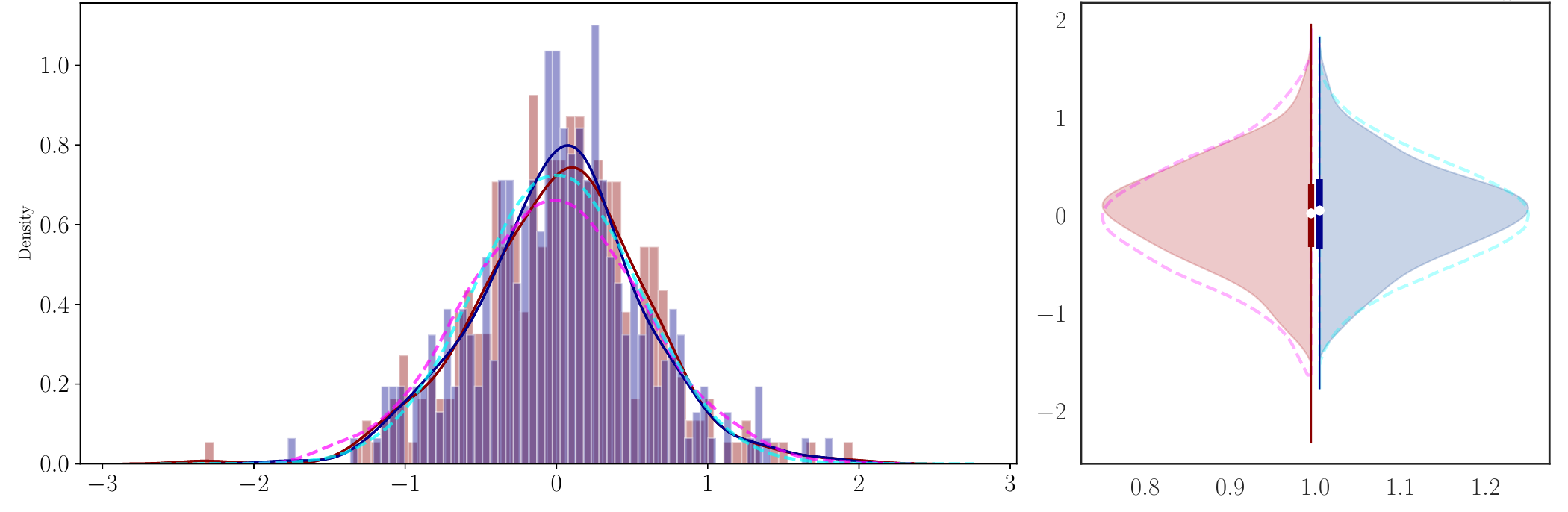}
    \end{minipage}
  
  \end{figure}
  \begin{figure}[H]
      \vspace{-0.8cm}
      \centering
      \begin{minipage}[b]{0.45\textwidth} % Adjust the width as needed
          \centering
          \includegraphics[width=7.2cm]{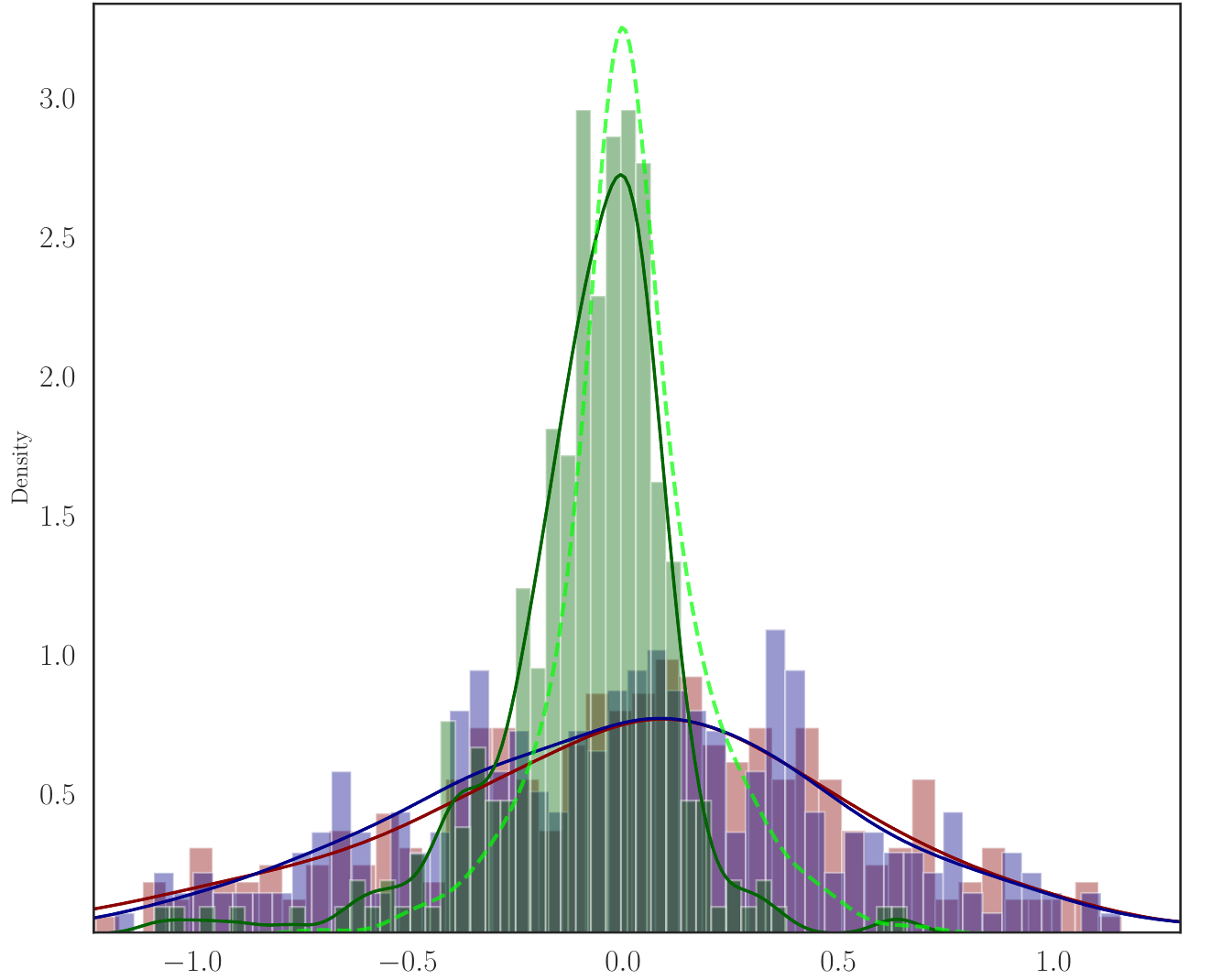}
          %\caption{First image}\label{fig:image1}
      \end{minipage}
      \hspace{0.8cm}
      \begin{minipage}[b]{0.45\textwidth} % Adjust the width as needed
          \centering
          \raisebox{2.5cm}{\includegraphics[width=4.2cm]{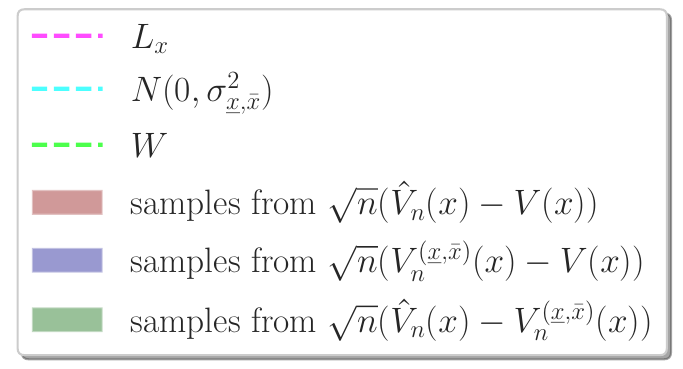}}
          %\caption{Second image}\label{fig:image2}
      \end{minipage}
      \caption{Based on $1000$ samples from $g$, centered and rescaled histograms based on $300$ repetitions, $V^{(\underline{x},\overline{x})}_n$ and their difference compared with respective limiting distributions $L_x$, $N(0,\sigma^2_{\underline{x},\overline{x}}) $ and $W$.}\label{img: histograms}
  \end{figure}
In Figure \ref{img: scatter + tuckey}, the scatter plot reveals a strong positive correlation between the two estimators, while the Tukey mean-difference plot indicates a lack of inherent bias. Taken together, the two plots in Figure \ref{img: scatter + tuckey} demonstrate a substantial agreement between the two estimation methods. 
  
In Figure \ref{img: histograms}, we observe that the distributions of $\sqrt{n} (V^{(\ubar{x},\bar{x})}_n(x)-V(x) )$ and $\sqrt{n} (\hat{V}_n(x)-V(x) )$ closely resemble each other. Moreover, we observe a strong resemblance not only between the estimators but also between their respective limiting distributions, $L_x$ and $N(0,\sigma^2_{\ubar{x},\bar{x}})$. However, the violin plot with the standard percentiles, the second histogram with $\sqrt{n}(\hat{V}_n(x) - V^{(\ubar{x},\bar{x})}_n(x))$ and $W$ and Table \ref{tab: 1} indicate that the variance of $\hat{V}_n(x)$ is slightly bigger than the one of $V^{(\ubar{x},\bar{x})}_n(x)$. 

The values in Table \ref{tab: 1} were computed using varying sample size $n$ from $g$ as indicated, with $2000$ samples from  $\hat{V}_n(x)$ and $V^{(\ubar{x},\bar{x})}_n(x)$. All the Kernel Density Estimator plots were computed using Scott's rule (where the bin size $h \sim 3.5 \hat{\sigma} n^{-1/3}$) for the chosen bandwidth. For $W$ and $L_x$, $20000$ samples were used.

Taken together, the simulation study confirms the theoretical findings of Theorem \ref{thm: distrib flat F at x total} and \ref{thm: LAM}: the informed estimators are efficient and not asymptotically equivalent to the IIE. However, the empirical performances of $V_n^{(\ubar{x},\bar{x})}$ and $\hat{V}_n$ show a high degree of proximity in the estimation scenario where the underlying function $V$ is constant on an interval. Even if $\hat{V}_n$ is used in cases where it is indeed true that $V$ is constant within an interval, the incurred loss appears to be small.

\begin{table}[H]
    \centering
    \resizebox{\columnwidth}{!}{%
    \begin{tabular}{|c|c c c c c||c|c|}
        \hline
         $n$ & 100 & 200 & 400 & 1000 & 2000  &  & limit \\
        \hline
         $\sqrt{n}(V^{(\ubar{x},\bar{x})}_n-V )(x)$ & 0.5441 & 0.5439 & 0.5465 & 0.5467 & 0.5464 & $N(0,\sigma^2_{\ubar{x},\bar{x}})$ & 0.5466 \\
        \hline
         $\sqrt{n} (\hat{V}_n-V)(x)$ & 0.5537 &  0.5592 & 0.5627 & 0.5664 & 0.5697  & $L_x$ & 0.5703 \\
        \hline
    \end{tabular}
    }
    \captionsetup{skip=6pt}
    \caption{Standard deviations}
    \label{tab: 1}
\end{table}

\vspace{-1cm}
\section{Conclusion}
We conclude by stressing the fact that the local asymptotically minimax variance in Theorem \ref{thm: LAM} coincides with the one obtained in Theorem \ref{thm: distrib flat F at x total 2}. This shows that in the setting in which $V$ is constant on $[\ubar{x},\bar{x}]$, for $\ubar{x} < \bar{x}$, the informed estimators \eqref{eq: empirical slope as proj}-\eqref{eq: proj isotonic smaller space} are asymptotically efficient, whereas the isotonic inverse estimator is not, as it is clearly shown by Theorem \ref{thm: LAM}. However, we conjecture that the isotonic estimator remains asymptotically efficient in the situations in which the boundary points $\ubar{x},\bar{x}$ are allowed to converge to $x$ at an arbitrarily polynomial rate (or slower). In this case, the rate of convergence drops to $\sqrt{n/\log{(\bar{x}-\ubar{x})^{-1}}}$. 

A potential direction for further investigation is whether the IIE is efficient in the setting in which there exists an interval of constancy for $V$ but the boundary points $\ubar{x},\bar{x}$ are unknown. However, the definition and assessment of \textit{efficiency} in this context are unclear. For practical purposes, if the boundary points $\ubar{x},\bar{x}$ are unknown, the IIE remains the state-of-art as it behaves very closely to the efficient estimator (see Section~\ref{sec:simulation study}) and in the $\sqrt{n / \log{n}}$ rate setting is efficient and adaptive to the level of smoothness of $V$.

%% The Appendices part is started with the command \appendix;
%% appendix sections are then done as normal sections
\newpage
\begin{appendix}

\section{Proofs.}
    
\section*{ Proofs of Section \ref{sec: construction estimators}.}\label{subsec: proofs lemmas thm normality holder}

To prove Proposition \ref{prop: construction estimators}, we need an additional lemma that will be essential to show that the estimators \eqref{eq: projection iso}-\eqref{eq: projection naive} are indeed the projections. The proof of this lemma uses similar reasoning as contained in \cite{23} and \cite{24}. Below for a function $f$, $f(t^{-})$ indicates the left limit of $f$ at $t$.
\begin{lemma}\label{lemma: minimization L^2}
Let $\overline{V} \in \mathbb{L}_1$ be given, define $V \mapsto Q^{\overline{V}}(V)$ by \eqref{eq: Q operator} and $\hat{V} \in \mathcal{V}_{\scriptstyle{\ubar{x},\bar{x}}}$.
Suppose the following conditions hold true, for $\ubar{x} < \bar{x}$: 
    \begin{enumerate}
        \item $ \int_{0}^{t} (\hat{V}(x) - \overline{V}(x) ) \, dx \geq 0$,  $\quad \forall \, t \in [0,\infty) \cap [\ubar{x},\bar{x}]^c$,
        \item $ \int_{0}^{t} ( \hat{V}(x) - \overline{V}(x) ) \, dx = 0$, $\quad \forall \, t \in [0,\infty) \cap [\ubar{x},\bar{x}]^c$ with $\hat{V}(t)<\hat{V}(t^{-})$.
        \item The measure defined via $\hat{V}$ and the one defined via $\overline{V}$ are compactly supported on the positive real line.
    \end{enumerate}
    Then $\hat{V} = \underset{V \in \mathcal{V}_{\scriptstyle{\ubar{x},\bar{x}}}}{\argmin} \: Q^{\overline{V}}(V)$.
\end{lemma}
\begin{proof}
Let us recall the concept of Gateaux differentiability at a point $V \in \mathcal{V}$, for some function space $\mathcal{V}$. Let $Q: \mathcal{V} \mapsto \mathbb{R}$ be an arbitrary functional. It is called Gateaux differentiable at the point $V \in \mathcal{V}$ if the limit:
$$
 \partial Q (V, h) =\lim _{\varepsilon \rightarrow 0} \frac{Q(V+\varepsilon h)-Q(V)}{\varepsilon}
$$
exists for every $h$ such that $V+\varepsilon h \in \mathcal{V}$ for small enough $\varepsilon$. By straightforward computations, one can verify that:
    \begin{align*}
        Q^{\overline{V}}(\hat{V} + \varepsilon \ind_{[0,t]}) - Q^{\overline{V}}(\hat{V}) = 2 \varepsilon \int_{0}^{t} \hat{V}(x) \, dx + \varepsilon^2 t - 2 \varepsilon \int_{0}^{t} \overline{V}(x) \, dx
    \end{align*}
    Therefore the Gateaux derivative of $Q^{\overline{V}}$ is given by:
    \begin{align}\label{eq: gateaux derivative}
        \partial Q^{\overline{V}} (\hat{V}, h) = 2 \int_{0}^{\infty} h(x) \left(\hat{V}(x) - \overline{V}(x) \right) \, dx
    \end{align}
    Note that, since $\partial Q^{\overline{V}} (\hat{V}, \ind_{[0,t]}) = \int_{0}^{t} \hat{V}(x) - \overline{V}(x) \, dx$, the first two conditions of Lemma \ref{lemma: minimization L^2} entail the Gateaux derivative of the functional $Q^{\overline{V}}$ along the test functions $\ind_{[0,t]}, \: \forall \, t \geq 0$ but $t \notin [\ubar{x},\bar{x}]$.  
    
    Now we show that for any $V \in \mathcal{V}_{\scriptstyle{\ubar{x},\bar{x}}}$, $\partial Q^{\overline{V}} (\hat{V}, V-\hat{V}) \geq 0$. First notice:
    \begin{align*}
        \hat{V}(x)-\overline{V}(x) = \frac{d(\hat{U}-\overline{U})}{dx}(x)
    \end{align*}
    Without loss of generality, we assume that the two supports of the measures defined via $\hat{V}$ and $\overline{V}$ are contained in $[0,M_1]$ and $[0,M_2]$. Defining $M := \max \{M_1,M_2\}$. Using integration by parts and assumption 3, $\hat{U}(0)-\overline{U}(0) = \hat{U}(M)-\overline{U}(M) = 0$:
    \begin{align*}
        \partial Q^{\overline{V}} (\hat{V}, V-\hat{V}) &= 2 \int_{0}^{\infty}  \left( V(x) - \hat{V}(x) \right) \left(\hat{V}(x) - \overline{V}(x) \right) \, dx \\
        &= 2 \int_{0}^{M}  \left( V(x) - \hat{V}(x) \right) d(\hat{U}-\overline{U}) (x) \\
        &=2 \underbrace{\left( V(x) - \hat{V}(x) \right) \left(\hat{U}(x)-\overline{U}(x)\right) \bigg|_{0}^{M}}_{=0} + \\
        &\quad \quad \quad - 2\int_{0}^{M} \left(\hat{U}(x)-\overline{U}(x)\right) d\left( V(x) - \hat{V}(x) \right) 
    \end{align*}
    By assumption 2, the support of the measure $- d \hat{V}$ is included in the set $\{x \: : \: \overline{U}(x) = \hat{U}(x) \}$ (cf.\ Lemma 1 in \cite{23}), therefore: $\int_{0}^{\infty} \left(\hat{U}(x)-\overline{U}(x)\right) d \hat{V}(x) =0 $. On the other hand, because $V \in \mathcal{V}_{\scriptstyle{\ubar{x},\bar{x}}}$ the negative measure defined based on $V$ does not have mass in $[\ubar{x},\bar{x}]$ and then, by using assumption 1, we obtain:
    \begin{align*}
        \int_{0}^{\infty} \left(\overline{U}(x) - \hat{U}(x)\right) d V(x) \geq 0
    \end{align*}
    which proves indeed $\partial Q^{\overline{V}} (\hat{V}, V-\hat{V})\geq 0$. Now define the convex function $t \mapsto u(t) = Q^{\overline{V}}(\hat{V}+ t (V- \hat{V}))$ where $t \in [0,1]$. Using \eqref{eq: gateaux derivative} and what we just proved, we have that $u^{\prime}(0) =\partial Q^{\overline{V}} (\hat{V}, V-\hat{V}) \geq 0$. Since $u$ is convex, the fact $u^{\prime}(0) \geq 0$ implies $u(1)\geq u(0)$, which indeed is $Q^{\overline{V}}(V) \geq Q^{\overline{V}}(\hat{V})$.
\end{proof}

\begin{proof}[\hypertarget{proof of prop constructions}{\textbf{Proof of Proposition}} \ref{prop: construction estimators}] The proof of \eqref{eq: isotonic estimator as proj} is provided in \cite{1} (c.f.\ Lemma 2 in \cite{1}). 

\vskip 0.1cm

\textit{Proof of} \eqref{eq: proj isotonic smaller space}. First, we show that $\hat{V}^{\Pi}_n \in \mathcal{V}_{\scriptstyle{\ubar{x},\bar{x}}}$.
    By definition, for $x \in [\ubar{x},\bar{x}]$, for $\ubar{x} < \bar{x}$:
    \begin{align*}
        \hat{V}^{\Pi}_n(x) = \frac{U^*_n(\bar{x}) - U^*_n(\ubar{x})}{\bar{x}-\ubar{x}}
    \end{align*}
    Let $Z^*_k$ for $k \in \{1,...,m-1 \}$ be the points where $U^*_n$ changes slope in $[\ubar{x},\bar{x}]$. Denote by $Z^*_0$ the closest point where $U^*_n$ changes slope to the left of $\ubar{x}$. Similarly, denote by $Z^*_m$ the closest point where $U^*_n$ changes slope to the right of $\bar{x}$. Now consider the fact that:
    \begin{align*}
        U_n^*(\ubar{x}) &= \frac{U_n(Z^*_1) - U_n(Z^*_0)}{Z^*_1 - Z^*_0} (\ubar{x}-Z^*_1) + U_n(Z^*_1)\\
        U_n^*(\bar{x}) &= \frac{U_n(Z^*_m) - U_n(Z^*_{m-1})}{Z^*_m - Z^*_{m-1}} (\bar{x}-Z^*_{m-1}) + U_n(Z^*_{m-1})
    \end{align*}
    By telescoping we have:
    \begin{align*}
        U_n(Z^*_1) - U_n(Z^*_{m-1}) = - \sum_{k=2}^{m-1} U_n(Z^*_k) - U_n(Z^*_{k-1})
    \end{align*}
    And therefore all together we have on $[\ubar{x},\bar{x}]$, for $\ubar{x} < \bar{x}$:
    \begin{align*}
        \hat{V}^{\Pi}_n(x) := &\sum_{k=2}^{m-1} \frac{Z^*_k - Z^*_{k-1}}{\bar{x}-\ubar{x}} \frac{U_n(Z^*_k) - U_n(Z^*_{k-1})}{Z^*_k - Z^*_{k-1}} + \frac{Z^*_1-\ubar{x}}{\bar{x}-\ubar{x}}\frac{U_n(Z^*_1)- U_n(Z^*_0)}{Z^*_1 - Z^*_0}\\
        &+ \frac{\bar{x}-Z^*_{m-1}}{\bar{x}-\ubar{x}}\frac{U_n(Z^*_m) - U_n(Z^*_{m-1})}{Z^*_m - Z^*_{m-1}} \numberthis \label{eq: construction iso proj}
    \end{align*}
    The fact that $\hat{V}^{\Pi}_n \in \mathcal{V}_{\scriptstyle{\ubar{x},\bar{x}}}$ is now clear, because the computations given above prove that $\hat{V}^{\Pi}_n$ on $[\ubar{x},\bar{x}]$ is just a weighted average of the values chosen by the isotonic estimator over $[\ubar{x},\bar{x}]$. This ensures both the monotonicity constraint and right-continuity, other than the fact that $\hat{V}^{\Pi}_n$ is constant on $[\ubar{x},\bar{x}]$. Let us now verify condition 1 of Lemma \ref{lemma: minimization L^2}. Clearly, for $t \leq \ubar{x}$, the condition is satisfied as it is zero. Take $t \geq \bar{x}$. The primitive of $\hat{V}$ is $U^*_n$. We verify the main condition. We split the integral into three bits. Clearly $\int_{0}^{\ubar{x}} \hat{V}^{\Pi}_n(x) - \hat{V}_n(x) \, dx = \int_{\bar{x}}^{t} \hat{V}^{\Pi}_n(x) - \hat{V}_n(x) \, dx = 0 $. Now for $[\ubar{x},\bar{x}]$: $\int_{\ubar{x}}^{\bar{x}} \hat{V}^{\Pi}_n(x) - \hat{V}_n(x) \, dx = 0$ by \eqref{eq: projection iso}. Condition 2 can be verified immediately using the same strategy and condition 3 is verified as well by the definition of the isotonic estimator.

    \vskip 0.1cm

\textit{Proof of} \eqref{eq: empirical slope as proj}. If we minimize $Q^{V_n}$ over a space where we do not require the solution to be decreasing (but still constant on each $[x_i,x_{i+1})$), then we can split the minimization problem into different integrals over each $[x_i,x_{i+1})$. When we minimize over $[x_i,x_{i+1})$ we obtain, setting $V(x)=c \in \mathbb{R}$ for $x \in [x_i,x_{i+1})$:
$$
\begin{aligned}
& \int_{x_i}^{x_{i+1}} V(x)\left(V(x)-2 V_n(x)\right) d x=\int_{x_i}^{x_{i+1}} c^2-2 c V_n(x) d x= \\
& =(x_{i+1}-x_i) c^2-2 c \frac{1}{n} \sum_{i=1}^n \int_{x_i}^{x_{i+1}} \left(Z_i-x\right)^{-1 / 2} d x \\
&=(x_{i+1}-x_i) c^2-2 c\left(U_n(x_{i+1})-U_n(x_i)\right).
\end{aligned}
$$
By optimizing, we obtain $c=\frac{U_n(x_{i+1})-U_n(x_i)}{x_{i+1}-x_i}$.
\vskip 0.1cm

    \textit{Proof of} \eqref{eq: proj naive smaller space}.  We use Lemma \ref{lemma: minimization L^2}. By \eqref{eq: projection naive}, $V^{\Pi}_n \in \mathcal{V}_{\scriptscriptstyle{\ubar{x},\bar{x}}}$. Moreover, we have that condition 1 is obtained because it requires that  $\forall \: t \in [\ubar{x},\bar{x}]^c \cap [0,\infty)$:
    \begin{align*}
        \int_{0}^{t} V^{\Pi}_n (x) - V_n(x) \, dx = \widetilde{U}^*_n (t) - U_n(t) \geq 0
    \end{align*}
    but by construction $\widetilde{U}_n^*(t)\geq \widetilde{U}_n(t) = U_n(t)$ for $t \notin [\ubar{x},\bar{x}]$. For condition 2, if $ V^{\Pi}_n(t) <  V^{\Pi}_n(t^{-})$ then $t$ is a support point of the $\widetilde{U}^*_n$ and we also know by construction that all the support points of $\widetilde{U}^*_n$ satisfy $\widetilde{U}^*_n(t) = U_n(t)$. Condition 3 is trivially satisfied.
\end{proof}

\section*{Proofs of Section \ref{sec: asymp distr}.}\label{sec: proofs 3}

All the proofs that follow make use of the explicit constructions of the estimators \eqref{eq: isotonic estimator as proj}-\eqref{eq: proj isotonic smaller space} given in Proposition \ref{prop: construction estimators}. We present the proof of Theorem \ref{thm: distrib flat F at x total 1} first, as part of the proof of Theorem \ref{thm: distrib flat F at x total 2} can be viewed as a special case of it.

To prove Theorem \ref{thm: distrib flat F at x total 2} and Theorem  \ref{thm: distrib flat F at x total 1} we need the following additional lemma.

\begin{lemma}\label{lemma: conv weak argmax}
Let $\ell^{\infty}(K)$ be endowed with the uniform norm and $\forall \, s \in K$:
\begin{align*}
    Z_n(s) := \frac{\sqrt{n}}{\bar{x}-\ubar{x}} \bigg( &(U_n(s)- U_n(\ubar{x})) (\bar{x} - s) - (U_n(\bar{x}) - U_n(s))(s - \ubar{x}) \numberthis \label{eq: Z_n}\\
    & - (U(s)- U(\ubar{x})) (\bar{x} - s) - (U(\bar{x}) - U(s))(s - \ubar{x})\bigg).
\end{align*}
Then
\begin{align}\label{eq: convergence ell^infty Z}
    Z_n \rightsquigarrow \mathbb{Z} \quad \text{in} \: \: \ell^{\infty}(K),
\end{align}    
where $\mathbb{Z}$ is the zero mean Gaussian Process given in \eqref{eq: gauss process} and has a version with continuous sample paths, with unique point of maximum.
\end{lemma}
\begin{proof}
Because
\begin{align*}
    &(U_n(s)- U_n(\ubar{x})) (\bar{x} - s) - (U_n(\bar{x}) - U_n(s))(s - \ubar{x}) \\
    &= \frac{2}{n} \sum_{i=1}^n \Bigg\{\left( \sqrt{(Z_i-\ubar{x})_{+}} - \sqrt{(Z_i-s)_{+}}\right)(\bar{x}-s) \\
    &- \left( \sqrt{(Z_i-s)_{+}} - \sqrt{(Z_i-\bar{x})_{+}}\right)(s-\ubar{x})  \Bigg\},
\end{align*}
the stated convergence is equivalent to saying that the class of functions $z \rightarrow ( \sqrt{(z-\ubar{x})_{+}} - \sqrt{(z-s)_{+}})(\bar{x}-s) - ( \sqrt{(z-s)_{+}} - \sqrt{(z-\bar{x})_{+}})(s-\ubar{x}), \: s \in K$ is Donsker. From Lemma 2.6.16 in \cite{5}, it follows that the class of functions $\{ z- s \: : \: s \geq 0 \}$ is VC, next by Lemma 2.6.18 (ii) in \cite{5} we see that the class $\{ (z -s)_{+} \: : \: s \geq 0 \}$ is VC. Again by Lemma 2.6.18 (vii), (iv) and (v) in \cite{5} we conclude that the class of functions $\{ \sqrt{(z-x)_{+}} - \sqrt{(z-s)_{+}} \: : \: s \geq 0 \}$ is VC. Similarly, we conclude that the class of functions 
$z \mapsto \sqrt{(z-s)_{+}} - \sqrt{(z-\bar{x})_{+}}$ is VC. The two above-mentioned classes of functions have respective envelopes: $0 \leq \sqrt{(z-x)_{+}} - \sqrt{(z-s)_{+}} \leq \sqrt{(z-x)_{+}}$ and $0 \leq \sqrt{(z-s)_{+}} - \sqrt{(z-\bar{x})_{+}} \leq \sqrt{z}$. Therefore these classes are $P$-Donsker for any $P$ with $\int_{0}^{\infty} z dP(z) < \infty$. This is the case by assumption \eqref{eq: finite first moment}. Again by Lemma 2.6.16 in \cite{5} we have that the classes of functions: $\{(\bar{x}-s) \: : \: s \in K \}$ and $\{(\ubar{x}-s) \: : \: s \in K\}$ are VC. They have clearly square integrable envelope as $K = [\ubar{x},\bar{x}]$ is compact. 
We therefore conclude that the uniform entropy integral for the class of functions of our interest is finite, using Lemma 7.21 (i) and (ii) from \cite{22}.

To prove that $\mathbb{Z}$ possesses a version with continuous sample paths, it suffices to show:
\begin{align}\label{eq: path distance}
    &\ex \left( \mathbb{Z} (s) - \mathbb{Z} (t) \right)^2   \lesssim |s - t| + |s - t|^2.
\end{align}
This follows from the fact that the square root is Hölder continuous of degree $1/2$ and the fact that:  
\begin{align*}
    &\frac{1}{\bar{x}-\ubar{x}} \left(\sqrt{(z-\ubar{x})_{+}} - \sqrt{(z-s)_{+}})(\bar{x}-s) - ( \sqrt{(z-s)_{+}} - \sqrt{(z-\bar{x})_{+}})(s-\ubar{x}) \right) \\
    &= \sqrt{(z-\ubar{x})_{+}} - \sqrt{(z-s)_{+}} + \frac{1}{\bar{x}-\ubar{x}} (\sqrt{(z-\ubar{x})_{+}} - \sqrt{(z-\bar{x})_{+}})(\ubar{x}-s)
\end{align*}
The covariance function of $\mathbb{Z}$ is given by, for $Z_x := Z -x$:
\begin{align*}
       &\mathrm{Cov} \left( \mathbb{Z}(t),\mathbb{Z}(s) \right)  \\
       &= 4\mathrm{Cov} \Bigg( \frac{(\sqrt{\smash[b]{(Z_{\ubar{x}})}_{+}} - \sqrt{(Z_t)_{+}})(\bar{x}-t) - (\sqrt{(Z_t)_{+}} - \sqrt{(Z_{\bar{x}})_{+}} )(t-\ubar{x})}{\bar{x} - \ubar{x}}, \\
        &\quad \quad \quad \quad \quad \quad  \frac{(\sqrt{\smash[b]{(Z_{\ubar{x}})}_{+}} - \sqrt{(Z_s)_{+}})(\bar{x}-s) - (\sqrt{(Z_s)_{+}} - \sqrt{(Z_{\bar{x}})_{+}} )(s-\ubar{x})}{\bar{x} - \ubar{x}} \Bigg).
\end{align*}
This coincides with the covariance function of the process in \eqref{eq: gauss process}.

By \eqref{eq: path distance}, we deduce that $\mathrm{Var}(\mathbb{Z}(t) - \mathbb{Z}(s)) \neq 0$ for $s \neq t$. Since $\mathbb{Z}$ is indexed by a $\sigma$-compact metric space, we obtain from Lemma 2.6 in \cite{12} that the location of the maximum of the sample paths of $\mathbb{Z}$ is a.s.\ unique. 
\end{proof}

\begin{proof}[\hypertarget{proof of thm 1 a}{\textbf{Proof of Theorem}} \ref{thm: distrib flat F at x total 1}]\label{proof: thm 1.a}
For $a >0$, the "switch relation" (see \cite{2}) is gives:
\begin{align}\label{eq: switch relation}
    & \hat{V}_n (x) \leq a \quad \Leftrightarrow \quad {\argmax_{s \geq 0}}^{-} \left\{ U_n(s) - as \right\} \leq x. 
\end{align}
Here the argmax is defined as: 
$${\argmax_{s \geq 0}}^{-} \left\{ U_n(s) - as \right\} := \inf \left\{s \geq 0: U_n(s)-a s \text { is maximal}\right\}.$$
Therefore for each fixed $a \in \mathbb{R}$ we can, for all $n$ sufficiently large, write:
\begin{align*}
    &\sqrt{n} \left( \hat{V}_n (x) -V^{(\ubar{x},\bar{x})}_n (x) \right) \leq a \quad \\
     \quad \Leftrightarrow \quad  &{\argmax_{s \geq 0}}^{-} \left\{ U_n(s) - \left(V^{(\ubar{x},\bar{x})}_n (x)+ a/\sqrt{n} \right) s \right\} \leq x \quad \quad \\
    \quad \Leftrightarrow \quad  &{\argmax_{s \geq 0}}^{-} \big\{ \sqrt{n} (U_n (s) - U_n(\ubar{x}) - U (s) + U(\ubar{x}) ) + \\
    &\quad \quad \quad \quad - \sqrt{n} ((U_n(\bar{x}) - U_n(\ubar{x}))/(\bar{x} - \ubar{x}) - V(x)) (s - \ubar{x}) + \\
    &\quad \quad \quad \quad + \sqrt{n}(U(s) - U(\ubar{x}) - V(x)(s-\ubar{x}) )  -  a s  \big\} \leq x. 
\end{align*}
Here we used that the location of a maximum of a function is equivariant under translations and invariant under multiplication by a positive number and addition of a constant. We obtain that the above is equivalent to:
\begin{align*}
    {\argmax_{s \geq 0}}^{-} \left\{ Z_n(s) - a s - \sqrt{n} h_x(s) \right\} \leq x,
\end{align*}
where:
\begin{align*}
    &Z_n(s) := \sqrt{n} (U_n (s) - U_n(\ubar{x}) - U (s) + U(\ubar{x}) )  \\
    &\quad \quad \quad \quad - \sqrt{n} ((U_n(\bar{x}) - U_n(\ubar{x}))/(\bar{x} - \ubar{x}) - V(x)) (s - \ubar{x}),
    \\ &h_x(s) := U(\ubar{x}) - U(s) - V(x)(\ubar{x}-s).
\end{align*}
One can easily show that $Z_n$ can be written as in \eqref{eq: Z_n}. Therefore by Lemma \ref{lemma: conv weak argmax} and the above computations, we conclude that the limiting behavior of: 
$$\widetilde{Z}_n (s) :=  Z_n(s) - a s - \sqrt{n} h_x(s), $$ 
is determined by the process $\mathbb{W}$ defined by, for $\ubar{x}< \bar{x}$:
$$
\mathbb{W}(s) = \begin{cases}
        \mathbb{Z}(s) -as, & \text{for } s \in [\ubar{x},\bar{x}]\\
        -\infty, & \text{for } s \notin [\ubar{x},\bar{x}]. 
        \end{cases} 
$$
This is because $h_{x}(s) > 0$ for all $s \in (0 ,\ubar{x}) \cup (\bar{x},\infty)$, and therefore multiplied by $-\sqrt{n}$ will go to $-\infty$, while $h_x(s) = 0, \: \: \forall \: s \in [\ubar{x},\bar{x}]$, and thus on that interval the limiting behavior is completely determined by $\mathbb{Z}(s) -as$. 

Let:
\begin{align*}
&\hat{s}_n := {\argmax_{s \geq 0}}^{-} \left\{  Z_n(s) - a s - \sqrt{n} h_x(s)  \right\} \numberthis \label{eq: def hat{t}_n },\\
&\hat{s} := \argmax_{s \geq 0} \left\{\mathbb{W}(s) \right\} = \argmax_{s \in [\ubar{x},\bar{x}]} \left\{\mathbb{Z}(s) - as \right\}. \numberthis \label{eq: def hat{t} }
\end{align*}
To prove convergence in distribution of $\hat{s}_n$ to $\hat{s}$, we use the Portmanteau Lemma as done in \cite{21} (cf.\ proof of Theorem 3) and show for every closed subset $F$: \\
$\limsup_{n \rightarrow \infty} \pb \left( \hat{s}_n \in F \right) \leq \pb \left( \hat{s} \in F \right)$. As in \cite{21} we choose a sequence $\varepsilon_n \downarrow 0$ such that $\sqrt{n} h_x(\ubar{x}-\varepsilon_n) \rightarrow \infty$ and $\sqrt{n} h_x(\bar{x}+\varepsilon_n) \rightarrow \infty$. For $K_n := [\ubar{x}-\varepsilon_n,\bar{x}+\varepsilon_n]$ and for $K := [\ubar{x},\bar{x}]$:
\begin{align*}
    &\limsup_{n \rightarrow \infty} \pb \left(\hat{s}_n \in F  \right) \leq \limsup_{n \rightarrow \infty} \pb \left(\hat{s}_n \in F \cap K_n \right) + \limsup_{n \rightarrow \infty} \pb \left(\hat{s}_n \in K^{c}_n  \right)  \\
    & \leq \underbrace{\limsup_{n \rightarrow \infty} \pb \left( \sup_{s \in F \cap K_n} \widetilde{Z}_n(s) \geq \sup_{s \in K_n} \widetilde{Z}_n(s) \right)}_{(1)} + \underbrace{\limsup_{n \rightarrow \infty} \pb \left(\hat{s}_n \in K^{c}_n  \right).}_{(2)} \numberthis \label{eq: two terms to show Portmanteau}
\end{align*}
We show that (2) in \eqref{eq: two terms to show Portmanteau} goes to zero. If $\: \forall \: s \in K^c$ we have $Z_n$ as in \eqref{eq: Z_n}, then:
    \begin{align}\label{eq: conv K^c_n}
        \sup_{s \in K^c_n} \widetilde{Z}_n(s) = \sup_{s \in K^c_n}  \big\{ Z_n(s) - a s - \sqrt{n} h_x(s) \big\} \stackrel{\pb}{\rightarrow} - \infty
    \end{align}
Note that:
    \begin{align*}
        &\sup_{s \in K^c_n} \left\{ Z_n(s) - \sqrt{n} h_x(s)-as \right\} \leq \sup_{s \in K^c} Z_n(s) - \inf_{\substack{s < \ubar{x}-\varepsilon_n \\ s > \bar{x}+\varepsilon_n}} \left\{ \sqrt{n} h_x(s) + as\right\} \\
        &\leq \sup_{s \in K^c} \{\sqrt{n}( U_n(s) - U_n(\ubar{x}) - U(s) + U(\ubar{x}))\} - \inf_{\substack{s < \ubar{x}-\varepsilon_n \\ s > \bar{x}+\varepsilon_n}} \left\{ \sqrt{n} h_x(s) + as\right\} \\
        &\quad \quad + \sup_{s \geq 0} \left\{-\sqrt{n}\left(U_n(\bar{x}) - U_n(\ubar{x}) - U(\bar{x}) + U(\ubar{x}) \right)\frac{s-\ubar{x}}{\bar{x}-\ubar{x}}\right\}  \\
        & = \sup_{s \in K^c} \{\sqrt{n} (U_n(s) - U_n(\ubar{x}) - U(s) + U(\ubar{x}))\}  - \inf_{\substack{s < \ubar{x}-\varepsilon_n \\ s > \bar{x}+\varepsilon_n}} \left\{ \sqrt{n} h_x(s) + as\right\} \\
        & \quad \quad + \sqrt{n}\frac{U_n(\bar{x}) - U_n(\ubar{x}) - U(\bar{x}) + U(\ubar{x})}{\bar{x}-\ubar{x}} \ubar{x}
    \end{align*}
    where the first and the last term on the right-hand side are $O_p(1)$ by the convergence given in Lemma 1 in \cite{21}. The rest of the proof is the same as in \cite{21}, and it shows that: $- \inf_{\substack{s < \ubar{x}-\varepsilon_n \\ s > \bar{x}+\varepsilon_n}} \left\{ \sqrt{n} h_x(s) + as\right\} \rightarrow - \infty.$
Using \eqref{eq: conv K^c_n}:
\begin{align*}
    \pb \left(\hat{s}_n \in K^{c}_n  \right) &\leq \pb \bigg( \sup_{s \in K^c_n} \widetilde{Z}_n(s) \geq \sup_{s \in K_n} \widetilde{Z}_n(s) \bigg) \\
    &\leq \pb \bigg( \sup_{s \in K^c_n} \widetilde{Z}_n(s) \geq \sup_{s \in K} \big\{ Z_n(s) -as \big\} \bigg) \rightarrow 0. \numberthis \label{eq: probabilities hat{t}_n in K^c}
\end{align*}
The last convergence is a consequence of \eqref{eq: conv K^c_n} and the convergence $Z_n \rightsquigarrow \mathbb{Z}$ on $\ell^{\infty} (K)$, combined with the continuous mapping theorem (which implies that the term on the right-hand side of the last inequality is $O_p(1)$).

Now we argue the behavior of (1) in \eqref{eq: two terms to show Portmanteau}. For that, we need an additional convergence in probability given in \eqref{eq: conv prob argmax}:
\begin{align}\label{eq: conv prob argmax}
    \sup_{s \in F \cap K_n} \left\{Z_n(s) - \sqrt{n} h_x(s) -as\right\} - \sup_{s \in F \cap K} \left\{Z_n(s) -as\right\} \stackrel{\pb}{\longrightarrow} 0.
\end{align}
First note:
\begin{align*}
    0 &\leq \sup_{s \in F \cap K_n} \left\{Z_n(s) - \sqrt{n} h_x(s) -as\right\} - \sup_{s \in F \cap K} \left\{Z_n(s) -as\right\} \\
    & \leq \sup_{s \in F \cap K_n} \left\{Z_n(s) -as\right\} - \sup_{s \in F \cap K} \left\{Z_n(s) -as\right\}.
\end{align*}
Define the process, for $\ubar{x}<\bar{x}$:
\begin{align*}
    \widetilde{Z}_n(s) := \begin{cases}{Z_n(s)-as,} & \text { if } s \in [\ubar{x},\bar{x}], \\ Z_n(\ubar{x})-a\ubar{x}, & \text { if } s \in [0,\ubar{x}], \\ Z_n(\bar{x})-a\bar{x}, & \text { if } s \in [\bar{x},\infty). \end{cases}
\end{align*}
Because: $\sup_{s \in F \cap K_n} \widetilde{Z}_n(s) = \sup_{s \in F \cap K}  \widetilde{Z}_n(s) = \sup_{s \in F \cap K} \left\{Z_n(s) -as\right\}$, it follows that
\begin{align*}
    &\quad \bigg| \sup_{s \in F \cap K_n} \left\{Z_n(s) -as\right\} - \sup_{s \in F \cap K} \left\{Z_n(s) -as\right\} \bigg| \\
    &= \bigg| \sup_{s \in F \cap K_n} \left\{Z_n(s) -as\right\} - \sup_{s \in F \cap K_n} \widetilde{Z}_n(s) \bigg| \leq \sup_{s \in F \cap K_n} \big| Z_n(s) - as - \widetilde{Z}_n(s) \big| \\
    & \leq \sup_{s \in F \cap K_n \cap [0,\ubar{x}]} \big| Z_n(s) - as - \widetilde{Z}_n(\ubar{x}) \big| \vee \sup_{s \in F \cap K_n \cap [\bar{x},\infty)} \big| Z_n(s) - as - \widetilde{Z}_n(\bar{x}) \big| \rightarrow 0,
\end{align*}
where we used the asymptotic equicontinuity of the $Z_n$.
Using that:
\begin{align*}
    \sup_{s \in F \cap K} \left\{Z_n(s) -as\right\} \rightsquigarrow \sup_{s \in F \cap K} \left\{\mathbb{Z}_x(s) -as\right\},
\end{align*}
we conclude using Theorem 2.7 (iv) from \cite{6} and \eqref{eq: conv prob argmax}:
\begin{align*}
    \sup_{s \in F \cap K_n} \left\{Z_n(s) - \sqrt{n} h_x(s) -as\right\} \rightsquigarrow  \sup_{s \in F \cap K} \left\{\mathbb{Z}_x(s) -as\right\}.
\end{align*}
Term (1) in \eqref{eq: two terms to show Portmanteau} can be upper bounded by:
\begin{align*}
    &\limsup \pb \left( \sup_{s \in F \cap K_n} \left\{Z_n(s) - \sqrt{n} h_x(s) -as\right\} \geq \sup_{s \in K} \left\{ Z_n(s) -as \right\} \right) \\
    &\leq \pb \left( \sup_{F \cap K} \left\{\mathbb{Z}_x(s) -as\right\} \geq \sup_{s \in K} \left\{\mathbb{Z}_x(s) -as\right\} \right) \leq \underbrace{\pb (\hat{s} \in K^c)}_{=0} + \pb \left( \hat{s} \in F \right), \numberthis \label{eq: first term portmanteau}
\end{align*}
where we used the above derivations, together with Lemma \ref{lemma: conv weak argmax} and Theorem 2.7 (v) from \cite{6}. Therefore by \eqref{eq: two terms to show Portmanteau} combined with \eqref{eq: probabilities hat{t}_n in K^c}, \eqref{eq: first term portmanteau} and the Portmanteau Lemma we obtain the desired convergence in distribution.
\end{proof}

To prove Theorem \ref{thm: distrib flat F at x total 2} we need the following additional lemma, which uses similar arguments as in \cite{25}.
\begin{lemma}\label{lemma: U^* - U = o_p}
If $U$ is linear with positive slope on $[\ubar{x},\bar{x}]$, for $\ubar{x}<\bar{x}$, and there does not exist a larger interval containing $[\ubar{x},\bar{x}]$ on which $U$ is linear, then:
    $$U^*_n(\bar{x})- U_n(\bar{x}) = o_p(1/\sqrt{n}) \quad \text{and} \quad U^*_n(\ubar{x})- U_n(\ubar{x}) = o_p(1/\sqrt{n})$$
\end{lemma}
\begin{proof} Because $U_n(x) = \frac{2}{n} \sum_{i=1}^{n} \left\{ \sqrt{(Z_i)_{+}} - \sqrt{(Z_i - x)_{+}} \right\}$, the same arguments as in the proof of Lemma \ref{lemma: conv weak argmax} give in $\ell^{\infty}(0,\infty)$:
\begin{align}\label{eq: weak conv tilde process}
    \sqrt{n} \left(U_n - U \right) \rightsquigarrow \widetilde{\mathbb{Z}} 
\end{align}
for a Gaussian Process $\widetilde{\mathbb{Z}}$. By assumption $U$ is linear on $[\ubar{x},\bar{x}]$. Let $l(t) = \frac{U(\bar{x})-U(\ubar{x})}{\bar{x}-\ubar{x}}(t-\ubar{x}) + U(\ubar{x})$  be the line which is equal to $U(t)$ on $[\ubar{x},\bar{x}]$. Recall that the upper script $*$ denotes the operation of taking the least concave majorant, for instance $U^*$ is the least concave majorant of $U$. Then consider:
\begin{align*}
    &\sqrt{n} \left( U^*_n(\ubar{x}) - U_n(\ubar{x}) \right) = \sqrt{n}  \left[ \left\{ U_n -U + U - l + l\right\}^*(\ubar{x}) - U_n(\ubar{x}) \right] \\
    &=  \left[ \sqrt{n} \left\{ U_n -U \right\} + \sqrt{n} \left\{U - l \right\} \right]^*(\ubar{x}) - \sqrt{n} \left(U_n(\ubar{x}) - U(\ubar{x}) \right) \numberthis \label{eq: expression last line},
\end{align*}
where we use that the least concave majorant of the sum of any function $h$ and a linear function $l$ is equal to the sum of the linear function $l$ and the least concave majorant of $h$. 

By Skorokhod's representation theorem there exists a sequence of random variables $\widetilde{Z}_n$, with same distribution as $\sqrt{n} \left\{ U_n -U \right\}$, which converges almost surely to $\widetilde{\mathbb{Z}}$ (or see \cite{34} for strong approximation type of result). Furthermore taking the Least Concave Majorant is a continuous function (see \cite{35}). Thus the expression in \eqref{eq: expression last line} has the same limiting distribution as:
\begin{align*}
    \left[ \widetilde{\mathbb{Z}} + \sqrt{n} \left\{U - l \right\} \right]^*(\ubar{x}) - \widetilde{\mathbb{Z}}(\ubar{x}).
\end{align*}
We proceed by an analytic argument to show that the above tends to zero in probability. Fix a sample path  $k$ of $\widetilde{Z}$, which is a bounded and continuous function on $[0,M]$. Because $0 \leq \sqrt{n}(U^*_n(\ubar{x}) - U_n(\ubar{x}))$, for $\varepsilon > 0$ we need to show $\forall \, n$ large enough:
\begin{align}\label{eq: preceding display}
    \{k+\sqrt{n}[U-l]\}^*(\ubar{x})<k(\ubar{x})+\varepsilon
\end{align}
To do so, we construct a concave majorant of $k + \sqrt{n}[U-l]$ on $(\ubar{x},M]$ passing through the point $(\ubar{x},k(\ubar{x}) + \varepsilon)$ and we show that for all $n$ big enough it constitutes a concave majorant of $k+ \sqrt{n}[U-l]$ also on $[0,\ubar{x}]$, thus eventually on the whole $[0,M]$.

Let $S_n$ denote the largest slope of lines connecting the point $(\ubar{x}, k(\ubar{x})+\varepsilon)$ with $(t, k(t)+\sqrt{n}[U-l](t))$ for $t>\ubar{x}$. Note that $S_n<\infty$ for all $n$. Since $U$ is concave, increasing and linear on $[\ubar{x},\bar{x}]$, $U-l \leq 0$ and thus $\sqrt{n}[U-l]$ becomes progressively smaller on $(\bar{x},M]$, implying $S_1 \geq S_2 \geq \cdots$.
 
Let $l_n$ be the line with slope $S_n$ that passes through the point $(\ubar{x}, k(\ubar{x})+\varepsilon)$. \eqref{eq: preceding display} holds true if and only if $k(t)+\sqrt{n}[U(t)-l(t)]<l_n(t)$ for all $t \in[0, \ubar{x}]$. Since $S_1 \geq S_n$, it suffices to show that for all $t \in[0, \ubar{x}]$ for all $n$ large enough
\begin{align}\label{eq: to prove l1}
k(t)+\sqrt{n}[U(t)-l(t)]<l_1(t),
\end{align}
as for $t > \ubar{x}$,  $k + \sqrt{n}[U-l] < l_1$ by construction and $S_1 \geq S_n$. These imply:
$$ \{ k + \sqrt{n}[U-l] \}^* \leq l_1.$$
We prove \eqref{eq: to prove l1}. Because $[U-l] \leq 0$, $\forall \, \delta >0$, as $n \rightarrow \infty$
$$ \sup_{0 \leq t < \ubar{x}-\delta} (k + \sqrt{n}[U-l])(t) \rightarrow - \infty$$
this implies that for any $\delta>0$, on $[0,\ubar{x}-\delta)$ for sufficiently large $n$:
$k + \sqrt{n}[U-l] < l_1$. For $\ubar{x}-\delta \leq t < \ubar{x}$ and $\delta < \ubar{x} - \varepsilon /(2 S_1)$, 
$$l_1(t) = k(\ubar{x}) + \varepsilon + (t- \ubar{x})S_1 \geq k(\ubar{x}) + \varepsilon/2.$$  
For $\delta$ small enough on $[\ubar{x}-\delta,\ubar{x})$ and sufficiently large $n$, by the continuity of $k$:
\begin{align*}
    \sup_{\ubar{x}-\delta \leq t < \ubar{x}} (k + \sqrt{n}[U-l])(t) \leq \sup_{\ubar{x}-\delta \leq t < \ubar{x}} k(t) \leq k(\ubar{x}) + \varepsilon/2.
\end{align*}
This implies that for $\delta$ small enough, on $[\ubar{x}-\delta,\ubar{x})$ for sufficiently large $n$: $k + \sqrt{n}[U-l] < l_1$. 

The result for $\bar{x}$ follows analogously.
\end{proof}

\begin{proof}[\hypertarget{proof of thm 1 b}{\textbf{Proof of Theorem}} \ref{thm: distrib flat F at x total 2}]
First, we show:
\begin{align}\label{first statement thm 1.b}
        \sqrt{n} \left(  V^{\Pi}_n(x) -V^{(\ubar{x},\bar{x})}_n (x) \right) = o_p(1),
\end{align}
Because $V_n^{\Pi}$ is the left derivative of the least concave majorant of $\widetilde{U}_n$, the switch relation gives, for $a \in \mathbb{R}$ and sufficiently large $n$,
\begin{align*}
    &\sqrt{n} \left(   V^{\Pi}_n(x) -V^{(\ubar{x},\bar{x})}_n (x) \right) \leq a \quad \\
    \quad \Leftrightarrow \quad  &{\argmax_{s \geq 0}}^{-} \big\{ \sqrt{n} (\widetilde{U}_n (s) - \widetilde{U}_n(\ubar{x}) - U (s) + U(\ubar{x}) ) + \\
    &\quad \quad \quad \quad - \sqrt{n} (V^{(\ubar{x},\bar{x})}_n (x) - V(x)) (s - \ubar{x}) + \\
    &\quad \quad \quad \quad + \sqrt{n}(U(s) - U(\ubar{x}) - V(x)(s-\ubar{x}) )  -  a s  \big\} \leq x,
\end{align*}
where we used the usual properties of the ${\argmax}^{-}$ as in the proof of Theorem \ref{thm: distrib flat F at x total 1}. By the definition of $\widetilde{U}_n$ in \eqref{eq: U_n tilde} and $V^{(\ubar{x},\bar{x})}_n (x)$ in \eqref{eq: natural competitor}, the above can be rewritten as:
\begin{align*}
    {\argmax_{s \geq 0}}^{-} \left\{ Z_n(s) - a s - \sqrt{n} h_x(s) \right\} \leq x,
\end{align*}
where:
\begin{align*}
    &Z_n(s) := \begin{cases}
        \frac{\sqrt{n}}{\bar{x}-\ubar{x}} \bigg( (U_n(s)- U_n(\ubar{x})) (\bar{x} - s) - (U_n(\bar{x}) - U_n(s))(s - \ubar{x}) \\
        \quad \quad - (U(s)- U(\ubar{x})) (\bar{x} - s) - (U(\bar{x}) - U(s))(s - \ubar{x})\bigg),   \quad \: s \notin [\ubar{x},\bar{x}],\\
        0,  \quad \quad \quad \quad \quad \quad \quad \quad \quad \quad \quad \quad \quad \quad \quad \quad \quad \quad \quad \quad \quad \quad  \quad \: \: s \in [\ubar{x},\bar{x}],
    \end{cases}\\ 
    &h_x(s) := U(\ubar{x}) - U(s) - V(x)(\ubar{x}-s).
\end{align*}
From here on, we apply the argument of Theorem \ref{thm: distrib flat F at x total 1}, using also for the current case the convergence in equation \ref{eq: conv K^c_n}. Presently
\begin{align*}
    {\argmax_{s \geq 0}}^{-} \left\{ Z_n(s) - a s - \sqrt{n} h_x(s) \right\} \rightsquigarrow \argmax_{s \in K} \left\{ -as \right\},
\end{align*}
as in this case the limiting process $\mathbb{Z}$ is constantly $0$ on $K$. The argmax on the right is $\bar{x}$ if $a<0$ and $\ubar{x}$ if $a>0$. Since $x \in (\ubar{x},\bar{x})$ it follows that
\begin{align*}
    \pb& \left(  {\argmax_{s \geq 0}}^{-} \left\{ Z_n(s) - a s - \sqrt{n} h_x(s) \right\} \leq x \right) \\
    &\stackrel{n \rightarrow \infty}{\longrightarrow}  \pb \left(  \argmax_{s \in K} \left\{ -as \right\} \leq x \right) = \begin{cases}
        1 \: \: \text{if} \: a>0, \\
        0 \: \: \text{if} \: a < 0.
    \end{cases}
\end{align*}
This concludes the proof of \eqref{first statement thm 1.b}. 

Second we prove:
\begin{align}\label{second statement thm 1.b}
        \sqrt{n} \left( \hat{V}^{\Pi}_n(x) -V^{(\ubar{x},\bar{x})}_n (x) \right) = o_p(1).
\end{align}
By the definitions of $\hat{V}^{\Pi}_n$ and $V^{(\ubar{x},\bar{x})}_n$ in \eqref{eq: projection iso} and \eqref{eq: natural competitor}:
    \begin{align*}
        \sqrt{n} &\left( \hat{V}^{\Pi}_n(x)-V^{(\ubar{x},\bar{x})}_n (x) \right) \hspace{-0.05cm} = \hspace{-0.05cm}\frac{1}{\bar{x}-\ubar{x}} \bigg\{ \hspace{-0.1cm} \sqrt{n} \left(  U^*_n(\bar{x}) - U_n(\bar{x}) \right) + \sqrt{n} \left(  U^*_n(\ubar{x}) - U_n(\ubar{x})  \right) \hspace{-0.1cm} \bigg\}.
    \end{align*}
Both terms on the right tend to zero by Lemma \ref{lemma: U^* - U = o_p}.

We conclude by giving the asymptotic distributions. By the Central Limit Theorem:
\begin{align}\label{eq: asymptotic distribution empirical slope}
       \sqrt{n} \left(V^{(\ubar{x},\bar{x})}_n (x)-V(x) \right) \rightsquigarrow N(0,\sigma^2_{\ubar{x},\bar{x}}).
\end{align}
The estimators $V^{\Pi}_n(x)$ and $\hat{V}^{\Pi}_n(x)$ have the same limiting behaviour by Slutsky's lemma and \eqref{first statement thm 1.b}-\eqref{second statement thm 1.b}.
\end{proof}

The \hyperlink{Proof of asymp indep}{proof} Theorem \ref{thm: asymptotic indep} is deferred to the next section because it relies on Lemma \ref{lemma: hadamard derivative}, Lemma \ref{lemma: IIE regular} and Proposition \ref{prop: LAN}, which can be understood only in the context of Section \ref{sec: local minimax}.

\section*{ Proofs of Section \ref{sec: local minimax}.}

\begin{proof}[\hypertarget{Proof of LAN prop}{\textbf{Proof of Proposition \ref{prop: LAN}}}]
    Define: 
    \begin{align*}
        &J_{n,i} := \left( 1 - \frac{1}{c_{t_n}} \right)+\frac{1}{g_{\scriptscriptstyle{V}}(Z_i)} \int_{\substack{\phantom{space} \\ \{   s \in [Z_i, M], \, (1+t_n h_0(s)) \leq 0\}}} \hspace{-2.5cm} c_{t_n}^{-1} \left(1+t_n h_0(s)\right) \frac{dV(s)}{\pi \sqrt{s-Z_i}}, \numberthis \label{eq: remainder ratio} \\
        &\Delta_{n,i} := \frac{1}{c_{t_n}\sqrt{n}} \left(k^{(\ubar{x},\bar{x})}(Z_i) - \mathbb{E}_{\scriptscriptstyle{G_{\scriptscriptstyle{V}}}} k^{(\ubar{x},\bar{x})}(Z) \right).
    \end{align*}
    Using the definition of $h_0$: 
    \begin{align*}
        &\frac{t_n}{c_{t_n}} \int_{Z_i}^{M} - \frac{h_0(s) v(s)}{\pi \sqrt{s-Z_i}} \, ds = t \Delta_{n,i} \, g_{\scriptscriptstyle{V}}(Z_i). \numberthis \label{eq: inversion formula for h_0} 
    \end{align*}
    Now note that for any $i$, using $1-\frac{(1+x)_{+}}{c} = - \frac{x}{c} + 1 - \frac{1}{c} + \frac{(1+x)_{-}}{c}$ and \eqref{eq: remainder ratio}-\eqref{eq: inversion formula for h_0}:
    \begin{align*}
        \frac{g_{\scriptscriptstyle{V_{t_n}}}(Z_i)}{g_{\scriptscriptstyle{V}}(Z_i)}  -1 &=\frac{1}{g_{\scriptscriptstyle{V}}(Z_i)}\left(\int_{Z_i}^M\left(1-\frac{\left(1+t_n h_0(s)\right)_{+}}{c_{t_n}}\right) \frac{d V(s)}{\pi \sqrt{s-Z_i}}\right) \numberthis \label{eq: second eq rewrite ratio}\\
        &=   - \frac{t_n}{c_{t_n} g_{\scriptscriptstyle{V}}(Z_i) }\int_{Z_i}^{M} - \frac{h_0(s) v(s)}{\pi \sqrt{s-Z_i}} \, ds +   J_{n,i} =  t \Delta_{n,i} +   J_{n,i}.
    \end{align*}
    Through the Taylor expansion $\log{\{1+x\}} = x - \frac{x^2}{2} + x^2 R(2x)$ where $R(x) \rightarrow 0$ as $x \rightarrow 0$, we can express the log-likelihood ratio:
    \begin{align*}
        &\sum_{i=1}^n \log{\frac{g_{\scriptscriptstyle{V_{t_n}}}(Z_i)}{g_{\scriptscriptstyle{V}}(Z_i)}} \numberthis \label{eq: expansion log} \\
        & = \sum_{i=1}^n \left\{ t \Delta_{n,i} +   J_{n,i} -\frac{1}{2}(t \Delta_{n,i} +   J_{n,i})^2 + (t \Delta_{n,i} +   J_{n,i})^2  R(2(t \Delta_{n,i} +   J_{n,i})) \right\} .
    \end{align*}
    We start by showing that:
    \begin{align}\label{eq: order of c_t_n}
        c_{t_n} =-\frac{2}{\pi} \int_0^M \sqrt{s}\left(1+t_n h_0(s)\right)_{+} d V(s) = 1 + o\left(\frac{1}{n} \right).
    \end{align}
    Let $S \sim \mu_{\scriptscriptstyle{V}}$ for $\mu_{\scriptscriptstyle{V}}$ defined in \eqref{eq: measure mu}. Since $\int h_0 \, d\mu_{\scriptscriptstyle{V}} = 0$:
    \begin{align*}
        c_{t_n} &  \geq -\frac{2}{\pi} \int_0^M \sqrt{s}\left(1+t_n h_0(s)\right) \, d V(s) =1.
    \end{align*}
     Now note:
\begin{align*}
c_{t_n} &  =-\frac{2}{\pi} \int_{\substack{\phantom{space} \\ \{   s \in [0, M], \, (1+t_n h_0(s)) > 0\}}} \hspace{-3cm} \sqrt{s}\left(1+t_n h_0(s)\right) d V(s)=1 + \frac{2}{\pi}  \int_{\substack{\phantom{space} \\ \{   s \in [0, M], \, (1+t_n h_0(s)) \leq 0\}}} \hspace{-3cm} \sqrt{s}\left(1+t_n h_0(s)\right) d V(s) \\
& =1+\mathbb{E}_{\scriptstyle{\mu_{\scaleto{V}{2.5pt}}}} \hspace{-0.1cm} \left[\left|1+t_n h_0(S)\right| \ind_{ \{  t_n h_0(S) \leq -1 \} } \right].
\end{align*}
If $t_n h_0(s) \leq -1$ then: $|1+ t_n h_0(s)| \leq |t_n h_0(s)| \leq t^2_n (h_0(s))^2$, so that:
\begin{align}\label{eq: upper bound h_0 on small set}
    \frac{1}{t_n^2} \ind_{ \{  t_n h_0(s) \leq -1 \} } \left|1+t_n h_0(s)\right| \leq (h_0(s))^2 \ind_{ \{  t_n h_0(s) \leq -1 \} }.
\end{align}
Now $\mathbb{E}_{\scriptstyle{\mu_{\scaleto{V}{2.5pt}}}} \hspace{-0.1cm} \left[ h^2_0(S)\right] \leq  \mathbb{E}_{\scriptstyle{\mu_{\scaleto{V}{2.5pt}}}} \hspace{-0.1cm} \left[ h^2(S)\right]$, which by \eqref{eq: bound on h} is bounded by \eqref{eq: L^2 condition}, which is finite by assumption. Therefore:
\begin{align}\label{eq: residual under measure mu_v}
    \limsup_{n \rightarrow \infty} \frac{1}{t^2_n} \mathbb{E}_{\scriptstyle{\mu_{\scaleto{V}{2.5pt}}}} \hspace{-0.1cm} \left[\left|1+t_n h_0(S)\right| \ind_{ \{  t_n h_0(S) \leq -1 \} } \right] = 0,
\end{align}
by the Dominated Convergence theorem. This concludes the proof of \eqref{eq: order of c_t_n} as $t^2_n = t^2/n$.

By \eqref{eq: order of c_t_n} and the Central Limit theorem: $\sum_{i=1}^n \Delta_{n,i} \rightsquigarrow N(0,\sigma^2_{\ubar{x},\bar{x}})$. Now note that because $V$ is a nonpositive measure, $\forall \, i$:
\begin{align}\label{eq: nonnegativity}
    \frac{1}{g_{\scriptscriptstyle{V}}(Z_i)} \int_{\substack{\phantom{space} \\ \{   s \in [Z_i, M], \, (1+t_n h_0(s)) \leq 0\}}} \hspace{-2.5cm} c_{t_n}^{-1} \left(1+t_n h_0(s)\right) \frac{dV(s)}{\pi \sqrt{s-Z_i}} \geq 0, \quad \quad \text{a.s.} 
\end{align}
Again by \eqref{eq: order of c_t_n}, by Markov's inequality and Fubini's theorem, $\sum_{i=1}^n J_{n,i} = o_p (1)$ because $\forall \, \varepsilon >0$:
\begin{align*}
&\mathbb{P}_{\scriptscriptstyle{G_{\scriptscriptstyle{V}}}}\left(\sum_{i=1}^n \frac{1}{g_{\scriptscriptstyle{V}}(Z_i)} \int_{\substack{\phantom{space} \\ \{   s \in [Z_i, M], \, (1+t_n h_0(s)) \leq 0\}}} \hspace{-2.5cm} c_{t_n}^{-1} \left(1+t_n h_0(s)\right) \frac{dV(s)}{\pi \sqrt{s-Z_i}} > \varepsilon \right) \\
&\leq \frac{n}{\varepsilon} \int_{0}^M \int_{z}^M \ind_{\{(1+t_n h_0(s)) \leq 0 \} } \, c_{t_n}^{-1} \left(1+t_n h_0(s)\right) \frac{v(s)}{\pi \sqrt{s-z}} \, ds \, dz \\
&=  \frac{n}{\varepsilon \pi} \int  \ind_{\{(1+t_n h_0(s)) \leq 0 \} } \, c_{t_n}^{-1} \left(1+t_n h_0(s)\right) v(s) \int_{0}^s \frac{1}{\sqrt{s-z}} \, dz \, ds \\
&= \frac{n}{\varepsilon c_{t_n}} \mathbb{E}_{\scriptstyle{\mu_{\scaleto{V}{2.5pt}}}} \hspace{-0.1cm} \left[\left|1+t_n h_0(S)\right| \ind_{ \{  t_n h_0(S) \leq -1 \} } \right] \stackrel{ \eqref{eq: residual under measure mu_v}}{=} o(1). \numberthis \label{eq: o(1) residual}
\end{align*}
All together this proves, for $\Delta_n$ as in \eqref{eq: linear expansion}: $\sum_{i=1}^n (t\Delta_{n,i} + J_{n,i}) = t\Delta_n + o_p(1)$, proving the behaviour of the linear term. 

For the quadratic term, if $\max_{1\leq i \leq n} |t \Delta_{n,i} + J_{n,i}| \stackrel{\mathbb{P}}{\rightarrow} 0 $, by property of the function $R$, the sequence 
$\max_{1\leq i \leq n} |R (t \Delta_{n,i} + J_{n,i})| \stackrel{\mathbb{P}}{\rightarrow} 0 $ as well. By the triangle inequality:
\begin{align*}
    \max_{1\leq i \leq n} |t \Delta_{n,i} + J_{n,i}| \leq |t| \max_{1\leq i \leq n}  |\Delta_{n,i}| + \sum_{i=1}^n |J_{n,i}|.
\end{align*}
The first term tends to zero in probability because for $\left\{ X_i , \, i \geq 1 \right\}$ i.i.d.\ random variables with finite variance: $\max_{1 \leq i \leq n} \left\{ \frac{|X_i|}{\sqrt{n}} \right\} \stackrel{\mathbb{P}}{\rightarrow} 0$. The second term tends to 0 by \eqref{eq: order of c_t_n} and \eqref{eq: o(1) residual}.

By the Law of Large numbers, $\sum_{i=1}^n \Delta^2_{n,i} \rightarrow \sigma^2_{\ubar{x},\bar{x}}$ in probability. Therefore we are left with showing:
\begin{align*}
    \sum_{i=1}^n \Delta_{n,i} J_{n,i} \stackrel{\mathbb{P}}{\rightarrow} 0 \quad \quad \text{and} \quad \quad  \sum_{i=1}^n (J_{n,i})^2 \stackrel{\mathbb{P}}{\rightarrow} 0. \numberthis \label{eq: squared o(1)} 
\end{align*}
The second follows because $ \sum_{i=1}^n (J_{n,i})^2 \leq \max_{1 \leq i \leq n} |J_{n,i}| \sum_{i=1}^n J_{n,i}$ and the first next follows with the help of the Cauchy-Schwarz inequality.
\end{proof}

For the proof of the LAM (locally asymptotically minimax) theorem, the first statement of Theorem \ref{thm: LAM}, and for the proof of Theorem \ref{thm: asymptotic indep}, we need the derivative of the map we are interested in estimating, in our case $V(x)$ (cf.\ definition 1.10 in \cite{13}). This is given in the following lemma.

\begin{lemma}\label{lemma: hadamard derivative}
    Along the path given in \eqref{eq: actual perturbation}, as $n \rightarrow \infty$, 
    \begin{align*}
        \sqrt{n} \left( V_{t_n}(x) - V(x)\right) \rightarrow  t \,  \sigma^2_{\ubar{x},\bar{x}}, 
    \end{align*}
    where $ \sigma^2_{\ubar{x},\bar{x}}$ is as in \eqref{eq: sigma}.
\end{lemma}
\begin{proof}
    Because $\int_{x}^s \frac{\sqrt{w-x}}{\sqrt{s-w}} \, dw = \frac{\pi}{2} (s - x)$ for every $x$, we have, $\forall \, x \in [\ubar{x},\bar{x}]$ by Fubini's theorem:
    \begin{align*}
        &\int_{0}^{M} 2 \frac{\sqrt{(w-\ubar{x})_{+}} - \sqrt{(w-\bar{x})_{+}}}{\bar{x}-\ubar{x}} \, g_{\scriptscriptstyle{V}} (w) \, dw \\
        &= - \frac{2}{\pi} \frac{1}{\bar{x}-\ubar{x}} \left( \int_{\ubar{x}}^{M}  \int_{\ubar{x}}^s \frac{\sqrt{w-\ubar{x}}}{\sqrt{s-w}} \, dw\, dV(s)  - \int_{\ubar{x}}^{M} \int_{\bar{x}}^s \frac{\sqrt{w-\bar{x}}}{\sqrt{s-w}} \, dw  \, dV(s) \right)\\
        &= - \frac{2}{\pi} \frac{1}{\bar{x}-\ubar{x}} \int_{\ubar{x}}^{M} \frac{\pi}{2} (\bar{x} - \ubar{x}) \, dV(s) = V(\ubar{x}) = V(x). \numberthis \label{eq: rewrite V(x)}
    \end{align*}
    Using the same computations as done for \eqref{eq: second eq rewrite ratio} we have:
    \begin{align*}
        & \left(g_{\scriptscriptstyle{V_t}}(z) - g_{\scriptscriptstyle{V}}(z) \right) =   - \frac{t_n}{ c_{t_n}} \int_{z}^{M} \frac{h_0(s) v(s)}{\pi \sqrt{s-z}} \, ds +  \left( 1 - \frac{1}{c_{t_n}} \right) \int_{z}^{M} \frac{dV(s)}{\pi \sqrt{s-z}} \\
        & \quad \quad \quad \quad \quad+  \int_{\substack{\phantom{space} \\ \{   s \in [z, M], \, (1+t_n h_0(s)) \leq 0\}}} \hspace{-2.5cm} c_{t_n}^{-1} \left(1+t_n h_0(s)\right) \frac{dV(s)}{\pi \sqrt{s-z}}.
    \end{align*}
    Because $(\mathbb{E}_{\scriptscriptstyle{G_{\scriptscriptstyle{V}}}} k^{(\ubar{x},\bar{x})} )\int_{0}^{M}(g_{\scriptscriptstyle{V_t}}(w) - g_{\scriptscriptstyle{V}}(w)) \, dw = 0$, using \eqref{eq: rewrite V(x)}, we obtain:
    \begin{align*}
        &\sqrt{n} \left( V_{t_n}(\ubar{x})- V(\ubar{x}) \right)= \sqrt{n} \int_{0}^{M} \left( k^{(\ubar{x},\bar{x})}(z) - \mathbb{E}_{\scriptscriptstyle{G_{\scriptscriptstyle{V}}}} k^{(\ubar{x},\bar{x})}  \right) \left( g_{\scriptscriptstyle{V_{t_n}}}(z) - g_{\scriptscriptstyle{V}}(z) \right) \, dz \\
        & =  \frac{t}{c_{t_n}} \int_{0}^{M} \left( k^{(\ubar{x},\bar{x})}(z) - \mathbb{E}_{\scriptscriptstyle{G_{\scriptscriptstyle{V}}}} k^{(\ubar{x},\bar{x})}  \right) \left(- \int_{z}^{M} \frac{h_0(s) v(s)}{\pi \sqrt{s-z}} \, ds \right) \, dz\\
        & \quad \quad + \sqrt{n}\left( \frac{1}{c_{t_n}} -1 \right) \underbrace{\int_{0}^{M} \left( k^{(\ubar{x},\bar{x})}(z) - \mathbb{E}_{\scriptscriptstyle{G_{\scriptscriptstyle{V}}}} k^{(\ubar{x},\bar{x})}  \right) g_{\scriptscriptstyle{V}} (z) \, dz}_{=0} \\
        & \quad \quad + \sqrt{n} \int_{0}^{M} \left( k^{(\ubar{x},\bar{x})}(z) - \mathbb{E}_{\scriptscriptstyle{G_{\scriptscriptstyle{V}}}} k^{(\ubar{x},\bar{x})}  \right) \int_{\substack{\phantom{space} \\ \{   s \in [z, M], \, (1+t_n h_0(s)) \leq 0\}}} \hspace{-2.5cm} c_{t_n}^{-1} \left(1+t_n h_0(s)\right) \frac{dV(s)}{\pi \sqrt{s-z}} \, dz\\
        &\stackrel{n \rightarrow \infty}{\longrightarrow} t \int_{0}^{M} \left( k^{(\ubar{x},\bar{x})}(w) - \mathbb{E}_{\scriptscriptstyle{G_{\scriptscriptstyle{V}}}} k^{(\ubar{x},\bar{x})}  \right)^2  g_{\scriptscriptstyle{V}}(w) \, dw = t \, \sigma^2_{\ubar{x},\bar{x}}.
    \end{align*}
   by \eqref{eq: inversion formula for h_0} and \eqref{eq: residual under measure mu_v}, where we use that $k^{(\ubar{x},\bar{x})}(z) - \mathbb{E}_{\scriptscriptstyle{G_{\scriptscriptstyle{V}}}} k^{(\ubar{x},\bar{x})}$ is bounded on $[0,M]$.
\end{proof}

The proof of Theorem \ref{thm: asymptotic indep}, is based on the fact that the isotonic estimator is regular on the submodel $t \mapsto G_{\scriptscriptstyle{V_t}}$, i.e.\ the limit distribution of the sequence $\sqrt{n}(\hat{V}_n (x) - V_{\scriptscriptstyle{G_{V_{t_n}}}} (x))$ under $G_{\scriptscriptstyle{V_{t_n}}}$ is the same, for all $t \in \mathbb{R}$.

\begin{lemma}\label{lemma: IIE regular}
    The IIE is regular in the submodel defined in \eqref{eq: actual perturbation}.
\end{lemma}
\begin{proof}
    For clarity of notation, in this proof we denote $U(x)$ as: $U_{\scriptscriptstyle{V}}(x) = \int_0^x V(y) \, dy$ and thus $U_{\scriptscriptstyle{V_{t_n}}}(x) = \int_0^x V_{t_n}(y) \, dy$. We show, under $G_{\scriptscriptstyle{V_{t_n}}}$: $\sqrt{n}(\hat{V}_n (x) - V_{t_n} (x)) \rightsquigarrow L_x$, which is independent of $t$. We use the same steps as in the proof of Theorem 4 in \cite{21}. Using the switch relation, we study: 
    \begin{align*}
        &\inf \big\{s \geq 0: (U_n (s) - U_n(x) - U_{\scriptscriptstyle{V_{t_n}}} (s) + U_{\scriptscriptstyle{V_{t_n}}}(x) ) + \\
    & \quad \quad \quad \quad + (U_{\scriptscriptstyle{V_{t_n}}}(s) - U_{\scriptscriptstyle{V_{t_n}}}(x) - V_{\scriptscriptstyle{t_n}}(x)(s-x) ) -  (a s)/\sqrt{n}  \text { is maximal} \big\} \leq x. 
    \end{align*}
    Using the same proof of Lemma 1 in \cite{21} it follows that, under $G_{\scriptscriptstyle{V_{t_n}}}$:
    \begin{align*}
        \sqrt{n}(U_n (s) - U_n(x) - U_{\scriptscriptstyle{V_{t_n}}} (s) + U_{\scriptscriptstyle{V_{t_n}}}(x) \, : \, s \geq 0) \rightsquigarrow \mathbb{Z}_x \quad \text{in} \: \: \ell^{\infty}[0,\infty).
    \end{align*}
    Furthermore:
    $$
    \begin{aligned}
        U_{\scriptscriptstyle{V_{t_n}}}(s)-U_{\scriptscriptstyle{V_{t_n}}}(x)-V_{t_n}(x)(s-x) & =\int_x^s\left(V_{t_n}(y)-V_{t_n}(x)\right) d y \\
        & = \begin{cases}0, & s, x \in[\ubar{x}, \bar{x}], \\
        <0, & x \leqslant \bar{x}<s, \quad s<\ubar{x} \leqslant x. \end{cases}
    \end{aligned}
    $$
    Now for $s = \bar{x} + \varepsilon_n$ and $\varepsilon_n \downarrow 0$ as in the proof of Theorem 4 in \cite{21} and using the same reasoning as in the proof of \eqref{eq: order of c_t_n}, we obtain:
    \begin{align*}
    &\sqrt{n}\int_{\bar{x}}^{\bar{x}+ \varepsilon_n}\left(V_{t_n}(y)-V_{t_n}(\bar{x})\right) d y \\
    & = -\sqrt{n} \int_{\bar{x}}^{\bar{x}+ \varepsilon_n} \int_{\bar{x}}^y c_{t_n}^{-1}\left(1+\frac{t}{\sqrt{n}} h_0(s)\right)_+ d V(s) d y \\
    & = -\sqrt{n} \int_{\bar{x}}^{\bar{x}+ \varepsilon_n} \int_{\bar{x}}^y c_{t_n}^{-1}\left(1+\frac{t}{\sqrt{n}} h_0(s)\right) d V(s) d y + o\left(\frac{1}{\sqrt{n}} \right)  \rightarrow - \infty,
    \end{align*}
    since by the proof of Theorem 4 in \cite{21} we have: $\sqrt{n}\int_{\bar{x}}^{\bar{x}+ \varepsilon_n} V(y) - V(\bar{x}) \, dy \rightarrow - \infty$. The same reasoning holds also for: $\sqrt{n} \int_{\ubar{x} - \varepsilon_n}^{\ubar{x}} \left( V_{t_n}(y)-V_{t_n}(\ubar{x})\right) \, dy \rightarrow - \infty$. Using further the same computations as in the proof of Theorem 4 in \cite{4} we conclude:
    \begin{align*}
        \pb \left(\sqrt{n} \big( \hat{V}_n (x) - V_{t_n} (x) \big) \leq a \right) \stackrel{n \rightarrow \infty}{\longrightarrow}  \pb \bigg( \argmax_{s \in [\ubar{x},\bar{x}]} \left\{ \mathbb{Z}_x(s) - a s\right\} \leq x \bigg).
    \end{align*}
\end{proof}

\vskip 0.1cm

\begin{proof}[\hypertarget{Proof of asymp indep}{\textbf{Proof of Theorem \ref{thm: asymptotic indep}}}](c.f.\ Theorem 1 in Section 2.3 of \cite{33}). 
    The efficient estimator satisfies:
    $$
    \sqrt{n}\big(V_n^{(\ubar{x}, \bar{x})}(x)-V_0(x)\big)=\Delta_n+o_{p}(1)
    $$
    where $V_0 \equiv V$ but this notation clarifies that it coincides with the perturbed $V_t$ in \eqref{eq: actual perturbation} evaluated at $t=0$. $\Delta_n \rightsquigarrow N(0, \sigma_{\ubar{x},\bar{ x}}^2)$ is the linear term in \eqref{eq: linear expansion} in the expansion of 
    $\Lambda_n (t) := \sum_{i=1}^n \log{\frac{g_{\scriptscriptstyle{V_{t_n}}}(Z_i)}{g_{\scriptscriptstyle{V}}(Z_i)}}$. The joint sequence
$$
\big(\sqrt{n} (\hat{V}_n(x)-V_n^{(x, x)}(x)), \sqrt{n}(V_n^{(x, x)}(x)-V_0(x))\big)
$$
is uniformly tight (because marginally so).
It suffices to show that every joint limit point under $V_0$ has independent coordinates. A joint limit point takes the form $(T-\Delta, \Delta)$ for $(T, \Delta)$ a joint limit of $\left(\sqrt{n}\left(\hat{V}_n(x)-V_0(x)\right), \Delta_n\right)$ under $V_0$. We have:
$$
\left(\sqrt{n}\left(\hat{V}_n(x)-V_0(x)\right), \Lambda_n(t)\right) \stackrel{V_0}{\rightsquigarrow}\bigg(T, t \Delta-\frac{1}{2} t^2 \sigma_{\ubar{x},\bar{x}}^2, \bigg).
$$
Therefore by Le Cam's 3rd lemma (e.g.\ Theorem 6.6 in \cite{6}):
$$
\begin{aligned}
& \sqrt{n}\left(\hat{V}_n(x)-V_0(x)\right) \stackrel{V_{t/\sqrt{n}}}{\rightsquigarrow} L_t, \quad \quad \quad L_t(B)=\mathbb{E} \ind_B(T) e^{t \Delta-\frac{1}{2} t^2 \sigma_{\ubar{x},\bar{x}}^2}
\end{aligned}
$$
because $\sqrt{n}\left(V_{t / \sqrt{n}}(x)-V_0(x)\right) \stackrel{n \rightarrow \infty}{\longrightarrow} t \sigma_{\ubar{x}, \bar{x}}^2$ (see Lemma \ref{lemma: hadamard derivative}) we conclude:
$$
\begin{aligned}
& \sqrt{n}\left(\hat{V}_n(x)-V_{t / \sqrt{n}}(x)\right)   \stackrel{V_{t/\sqrt{n}}}{\rightsquigarrow} \bar{L}_t, \quad \quad \quad \bar{L}_t(B)=L_t\left(B-t \sigma_{\ubar{x},\bar{x}}^2\right).
\end{aligned}
$$
By regularity of the IIE (see Lemma \ref{lemma: IIE regular}): $\bar{L}_t=\bar{L}_0$, independent of $t$, so $\forall \, t, u$
$$
\int e^{i u y} d \bar{L}_0(y)=\mathbb{E} e^{i u\left(T-t \sigma_{\ubar{x},\bar{x}}^2\right)} e^{t \Delta-\frac{1}{2} t^2 \sigma_{\ubar{x},\bar{x}}^2}.
$$
Choose $t=-i u$ to find: $\int e^{i u y} d \bar{L}_0(y)=\mathbb{E} e^{i u(T-\Delta)} e^{-\frac{1}{2} u^2 \sigma_{\ubar{x},\bar{x}}^2}$. Choose $t=-i u+i v$ to find: $\int e^{i u y} d \bar{L}_0(y)=\mathbb{E} e^{i u(T-\Delta)+i v \Delta} e^{\frac{1}{2}\left(v^2-u^2\right) \sigma_{\ubar{x}, \bar{x}}^2}$. Then:
$$
\begin{aligned}
\mathbb{E} e^{i u(T-\Delta)+i v \Delta}=\mathbb{E} e^{i u(T-\Delta)} e^{-\frac{1}{2} v^2 \sigma_{\ubar{x}, \bar{x}}^2} = \ex e^{i u(T-\Delta)} \ex e^{i v \Delta}. 
\end{aligned}
$$    
\end{proof}

\begin{proof}[\hypertarget{Proof of LAM thm}{\textbf{Proof of Theorem \ref{thm: LAM}}}]
The proof is a direct application of Theorems 25.20 and 25.21 from \cite{6}, using Proposition \ref{prop: LAN} and Lemma \ref{lemma: hadamard derivative}. The variance of the limiting random variable in \eqref{eq: LAM} is given by $\dot{\nu} (\sigma^2_{\ubar{x},\bar{x}})^{-1} \dot{\nu} = \sigma^2_{\ubar{x},\bar{x}}$, where $\dot{\nu} =  \sigma^2_{\ubar{x},\bar{x}}$ is the derivative of the functional of interest and $(\sigma^2_{\ubar{x},\bar{x}})^{-1}$ is the Fisher information for $t$, given in Proposition \ref{prop: LAN}.
\end{proof}

\section{Efficiency Theory}

In this section we derive the least favourable submodel \eqref{eq: actual perturbation} from the semiparametric score calculus, as in \cite{36} (or Chapter 25 of \cite{6}).

The parameter $V$ belongs to the class (c.f.\ \cite{1}), for $\ubar{x} < \bar{x}$:
\begin{align*}
&\mathscr{V}_{\ubar{x},\bar{x}}= \bigg\{ V \in \mathcal{V}_{\ubar{x},\bar{x}}, \: : \: \int_0^{\infty} \sqrt{s} \, d V(s)=-\frac{\pi}{2}, \:  \text{constant on } [M,\infty) \bigg\}. \numberthis \label{eq: def V space tangent}
\end{align*}
Because $V$ is not a probability measure, it is convenient to work with its associated probability measure $\mu_{\scriptscriptstyle{V}}$ defined for every measurable $A \subseteq [0,\infty)$ by:
\begin{align}\label{eq: measure mu}
    \mu_{\scriptscriptstyle{V}}(A) := \int_A - \frac{2}{\pi} \sqrt{s} \, dV(s).
\end{align}
We consider the associated space of "underlying" probability measures $\mathscr{M}_{\ubar{x},\bar{x}} := \{ \mu_{\scriptscriptstyle{V}} \, : \, V \in \mathscr{V}_{\ubar{x},\bar{x}} \}$ and the space of the "observed" probability measures $\mathscr{G}_{\ubar{x},\bar{x}} := \{ g_{\scriptscriptstyle{V}} \, : \, \mu_{\scriptscriptstyle{V}} \in \mathscr{M}_{\ubar{x},\bar{x}}\}$. The tangent set $\dot{\mathscr{M}}_{\ubar{x},\bar{x}}$ consists of all functions $h \in \mathbb{L}_2(\mu_{\scriptscriptstyle{V}})$ such that $\int h \, d \mu_{\scriptscriptstyle{V}} = 0$. The derivative of the map $\mu_{\scriptscriptstyle{V}} \mapsto g_{\scriptscriptstyle{V}}$ (the score operator) maps score functions in the tangent set $\dot{\mathscr{M}}_{\ubar{x},\bar{x}} $ into scores in the tangent set $\dot{\mathscr{G}}_{\ubar{x},\bar{x}} $. For any $h \in \mathbb{L}_1(\mu_{\scriptscriptstyle{V}})$ such that $- \int_{z}^{M} \frac{|h(s)|}{\sqrt{s-z}} \, dV(s) < \infty$, define:
\begin{align*}
    (Th)(z) := \int_{z}^{M} \frac{h(s)}{\sqrt{s-z}} \, dV(s). \numberthis \label{eq: operator T}
\end{align*}
Then the score operator $A_{\scriptscriptstyle{V}} \: : \: \dot{\mathscr{M}}_{\ubar{x},\bar{x}} \rightarrow \dot{\mathscr{G}}_{\ubar{x},\bar{x}}$, for any $h \in \dot{\mathscr{M}}_{\ubar{x},\bar{x}}$, is given by:
\begin{align}\label{eq: operator A_V}
    A_{\scriptscriptstyle{V}} h(z) = \frac{\int_z^{M} \frac{h(s)}{\sqrt{s-z}} \, d V(s)}{\int_z^{M} \frac{d V(s)}{\sqrt{s-z}}} = \frac{(Th)(z)}{(T 1)(z)}.
\end{align}
For instance, for bounded $h \in \dot{\mathscr{M}}_{\ubar{x},\bar{x}}$ define the following perturbation:
\begin{align}\label{eq: initial perturbation}
    dV_t(s) = (1+ t h(s)) \, dV(s).
\end{align}
Then:
$$
\begin{aligned}
& \left.\frac{d}{d t} \log \left\{- \int_z^{M} \frac{(1+ t h(s))}{\pi \sqrt{s-z}} d V(s)\right\}\right|_{t=0}=\frac{\int_z^{M} \frac{h(s)}{\sqrt{s-z}} d V(s)}{\int_z^{M} \frac{d V(s)}{\sqrt{s-z}}}=: A_{\scriptscriptstyle{V}} h(z). 
\end{aligned}
$$ 
We check that $A_{\scriptscriptstyle{V}}$ is a linear operator between the $\mathbb{L}_2$ spaces claimed above. Linearity is clear. We can check by using the Cauchy-Schwarz inequality and that $h \in \mathbb{L}^2(\mu_{\scriptscriptstyle{V}})$:
\begin{align*}
    &\int (A_{\scriptscriptstyle{V}} h)^2 dG_{\scriptscriptstyle{V}} = \int_0^{\infty}\left(\frac{\int_z^{\infty} \frac{h(s) d V(s)}{\sqrt{s-z}}}{\int_z^{\infty} \frac{d V(s)}{\sqrt{s-z}}}\right)^2 \left(-\frac{1}{\pi} \int_z^{\infty} \frac{d V(s)}{\sqrt{s-z}} \right) dz \numberthis \label{eq: score operator in L^2}\\
    &\leqslant \int_0^{\infty}\left(\frac{\sqrt{\int_{z}^{\infty} -\frac{h^2(s)}{\sqrt{s-z}} d V(s)} \sqrt{\int_z^{\infty}- \frac{1}{\sqrt{s-z}} d V(s)}}{-\int_z^{\infty} \frac{d V(s)}{\sqrt{s-z}}}\right)^2\left(-\frac{1}{\pi} \int_z^{\infty} \frac{d V(s)}{\sqrt{s-z}}  \right) dz \\
    &=-\frac{1}{\pi} \int_0^{\infty} \int_z^{\infty} \frac{h^2(s)}{\sqrt{s-z}} d V(s) d z=-\frac{2}{\pi} \int_0^{\infty} h^2(s) \sqrt{s} d V(s) = \int h^2 \, d\mu_{\scriptscriptstyle{V}}<\infty.
\end{align*}

The function $k^{(\ubar{x},\bar{x})}$ is the influence function of the (regular) estimator \eqref{eq: natural competitor}. To show that this estimator is asymptotically efficient, it is suffices to show that $k^{(\ubar{x},\bar{x})}$ is contained in the closed linear span of the set of all score functions $A_{\scriptscriptstyle{V}}h$. We shall show that in fact $k^{(\ubar{x},\bar{x})} - \ex k^{(\ubar{x},\bar{x})} = A_{\scriptscriptstyle{V}}h_0$ for $h_0$ given in \eqref{eq: h_0}. Equivalently, we show that $h_0$ solves
\begin{align}\label{eq: integral eq}
    \begin{cases}
        (Th_0) = \left( k^{(\ubar{x},\bar{x})} - \mathbb{E}_{\scriptscriptstyle{G_{\scriptscriptstyle{V}}}} k^{(\ubar{x},\bar{x})} \right) \cdot (T1) \\
        \int_{0}^{M}- | h_0(s)| \, dV(s) < \infty, \quad \quad \int_{0}^{M} \sqrt{s} h_0(s) \, dV(s)  = 0.
    \end{cases}
\end{align}
\begin{proposition}\label{prop: solution integral eq}
    Under Assumption \ref{assumption: LAM}, the solution to \eqref{eq: integral eq} is given by \eqref{eq: h_0}.
\end{proposition}

\begin{proof}[\hypertarget{proof of solution integral eq}{\textbf{Proof of Proposition}}  \ref{prop: solution integral eq}]
    We show that the function $h$ in \eqref{eq: h funtion} solves: $Th = q$, for $q = k^{(\ubar{x},\bar{x})} \cdot T1$. Then the solution $h_0$ is found by centering $h$, i.e.\ $h_0 (z) = h(z) - m$ (since $\int_{0}^{M} \sqrt{s} \, dV(s) = - \frac{\pi}{2}$) where $m$ is such that:
    \begin{align*}
        \int_{0}^{M} \sqrt{s} \left( h(s) - m \right) \, dV(s) = 0 \quad \implies \quad m = - \frac{2}{\pi} \int_{0}^M \sqrt{s} h(s) \, dV(s)
    \end{align*}
    The equation $Th = q$, for $q = k^{(\ubar{x},\bar{x})} \cdot T1$ is the following Abel's integral equation:
    \begin{align*}
        \int_{z}^{M} \frac{h(s)v(s)}{\sqrt{s-z}} \, ds = 2 \frac{\sqrt{(z-\ubar{x})_{+}} - \sqrt{(z-\bar{x})_{+}}}{\bar{x}-\ubar{x}} \int_{z}^{M} \frac{v(s)}{\sqrt{s-z}} \, ds
    \end{align*}
    By Theorem 2.1 in \cite{26} this equation is solvable for $hv \in \mathbb{L}_1[0,M]$ iff:
    \begin{align*}
        z \mapsto \int_{z}^{M} \frac{q(s)}{\sqrt{s-z}} \, ds
    \end{align*}
    is absolutely continuous. Furthermore the unique in $\mathbb{L}_1$ solution is given by $h(z) v(z) = - \frac{d}{dz} \int_{z}^{M} \frac{q(s)}{\sqrt{s-z}} \, ds$, where the derivative is understood in the sense of absolute continuity (AC). Moreover, if $q \in \operatorname{AC}[0,M]$, then this condition is satisfied and for $z \notin [\ubar{x},\bar{x}]$:
    \begin{align}\label{eq: hv rewriting}
        h(z)v(z) = \frac{1}{\sqrt{\pi}}\left( \frac{q(M)}{\sqrt{M-z}} - \int_{z}^{M} \frac{q^{\prime} (s)}{\sqrt{s-z}} \, ds \right).
    \end{align}
    Now we prove that under the given conditions $ k^{(\ubar{x},\bar{x})} \cdot T1 = q \in \operatorname{AC}[0,M]$. Because $\sqrt{(z-x)_{+}} = \int_{0}^{z} \frac{1}{2\sqrt{s-x}} \ind_{s>x} \, ds$ we see that $ k^{(\ubar{x},\bar{x})} \in \operatorname{AC}[0,M]$. $T1$ is the $1/2$-fractional integral of $v$. Because $v \in \operatorname{AC}[\ubar{x},\bar{x}]^c$ and $v=0$ on $[\ubar{x},\bar{x}]$, we know that $T1 \in \operatorname{AC}[0,M]$ (c.f.\ property 3.2.(6).(b) pp.\ 209 in \cite{27}). The product of absolutely continuous functions is absolutely continuous. 
    
    Finally we show:
    \begin{align*}
        \mathbb{E}_{\scriptscriptstyle{G_{\scriptscriptstyle{V}}}} k^{(\ubar{x},\bar{x})} = -\frac{2}{\pi} \int_{0}^M \sqrt{s} h(s) v(s) \, ds.
    \end{align*}
    For the given $h$ in \eqref{eq: h funtion} we just proved:
    \begin{align}\label{eq: mapping of h}
        \int_{z}^{M}  -\frac{h(s)v(s)}{\pi \sqrt{s-z}} \, ds = 2 \frac{\sqrt{(z-\ubar{x})_{+}} - \sqrt{(z-\bar{x})_{+}}}{\bar{x}-\ubar{x}} g_{\scriptscriptstyle{V}}(z).
    \end{align}
    By integrating both sides between $0$ and $M$ we get that the r.h.s.\ is $ \mathbb{E}_{\scriptscriptstyle{G_{\scriptscriptstyle{V}}}} k^{(\ubar{x},\bar{x})}$, whereas the l.h.s.\
    \begin{align*}
       \int_{0}^{M} \int_{z}^{M}  -\frac{h(s)v(s)}{\pi \sqrt{s-z}} \, ds \, dz = \int_{0}^{M} \int_{0}^{s}  -\frac{h(s)v(s)}{\pi \sqrt{s-z}} \, dz \, ds= -\frac{2}{\pi} \int_{0}^M \hspace{-0.2cm} \sqrt{s} h(s) v(s) \, ds.
    \end{align*}
    All together this proves that $h_0$ solves \eqref{eq: integral eq}.
\end{proof}

\begin{lemma}
    For a fixed $\eta >0$, there exist constants $c_1,c_2,c_3 >0$ such that, for $z \in [\ubar{x},\bar{x}]^c$:
    \begin{align}\label{eq: bound on h}
       |h(z)v(z)| \leq c_1  + c_2 \log{\left(\frac{1}{\ubar{x}-z}\right)} \ind_{\{ 0<\ubar{x}-z \leq \eta \}} + c_3 \log{\left(\frac{1}{z-\bar{x}}\right)} \ind_{\{0<z-\bar{x} \leq \eta \}}.
    \end{align}
\end{lemma}
\begin{proof}
Because $\lim_{z \rightarrow M}(T 1)(z) =0 $, the function $q$ in \eqref{eq: hv rewriting} satisfies $q(M) =0$. Hence, for $z \in [\ubar{x},\bar{x}]^c$:
\begin{align*}
  &h(z)  = -\frac{1}{\sqrt{\pi}v(z)} \Bigg( \underbrace{\int_{z}^M \frac{(k^{(\ubar{x},\bar{x})})^{\prime} (T1) (s)}{ \sqrt{s-z}} \, ds}_{(1)} + \underbrace{\int_{z}^M \frac{(k^{(\ubar{x},\bar{x})}) (T1)^{\prime} (s)}{ \sqrt{s-z}} \, ds }_{(2)} \Bigg) 
\end{align*}
where:
\begin{align*}
  \quad (k^{(\ubar{x},\bar{x})})^{\prime}(s) = \frac{1}{(\bar{x}-\ubar{x})} \left( \frac{1}{\sqrt{(s-\ubar{x})_{+}}} - \frac{1}{\sqrt{(s-\bar{x})_{+}}}\right).
\end{align*}

We prove (2) is bounded. By the same proof of Corollary 2.1 in \cite{26} (pp.\ 32), since $(T1) \in \operatorname{AC}[0,M]$ (see proof of Proposition \ref{prop: solution integral eq}), $\lim_{z \rightarrow M} v (z) =0 $, and denoting by $\operatorname{I}^{\frac{1}{2}}(f)$ the half--integral of $f$ (see \cite{26} for an introduction to fractional calculus):
$$\frac{d}{dz} \int_{z}^{M} \frac{v(s) \, ds}{\sqrt{s-z}} = \int_{z}^{M} \frac{v^{\prime}(s) \, ds}{\sqrt{s-z}} -  \lim_{s \rightarrow M}\frac{v(s) }{\sqrt{s-z}} = \operatorname{I}^{\frac{1}{2}}(v^{\prime}),$$
where the derivative is understood in the sense of absolute continuity (AC). Thus (2) coincides with $\operatorname{I}^{\frac{1}{2}}(k^{(\ubar{x},\bar{x})} \operatorname{I}^{\frac{1}{2}}(v^{\prime}))$. By Theorem 3.6 in \cite{26}, $\operatorname{I}^{\frac{1}{2}} : \mathbb{L}_p \mapsto H^{1/2 - 1/p}$, where $H^{\alpha}$ denotes the space of Hölder continuous functions of degree $\alpha$. Thus if $k^{(\ubar{x},\bar{x})} \operatorname{I}^{\frac{1}{2}}(v^{\prime}) \in \mathbb{L}_p$ for $p>2$ we obtain the claim. Because $k^{(\ubar{x},\bar{x})}$ is bounded we show $\operatorname{I}^{\frac{1}{2}}(v^{\prime}) \in \mathbb{L}_p$ for $p>2$. By theorem 3.5 in \cite{26},  $\operatorname{I}^{\frac{1}{2}}(v^{\prime}) \in \mathbb{L}_{q/(1-q/2)}$ if $v^{\prime} \in \mathbb{L}_q$. Because $v$ is Lipschitz on $[\ubar{x},\bar{x}]^c$ and $0$ on $[\ubar{x},\bar{x}]$, $ \exists \, \, q >1$ such that $v^{\prime} \in \mathbb{L}_q$. As $q>1 \iff q/(1-q/2) >2$, $\operatorname{I}^{\frac{1}{2}}(v^{\prime}) \in \mathbb{L}_p$ for some $p>2$ and we obtain the claim.

Regarding the absolute value of (1), using the fact that $(T1)(s) = -\frac{\pi}{2}g_{\scriptscriptstyle{V}}(s)$, which is bounded as well as $v$:
\begin{align*}
    \bigg| -\frac{1}{\sqrt{\pi}}  \int_{z \vee \ubar{x}}^M \frac{(T1) (s)}{ \sqrt{ (s-\ubar{x})(s-z)}} \, ds \bigg| \lesssim  \bigg| \int_{z \vee \ubar{x}}^M \frac{ds}{ \sqrt{ (s-\ubar{x})(s-z)}}   \bigg|.
\end{align*}
For $z > \bar{x}$ the previous display is bounded. For $z < \ubar{x}$ it is upper bounded (up to constants) by:
\begin{align*}
    \bigg| \sin^{-1} \left( \sqrt{\frac{M-z}{\ubar{x}-z}} \right) \bigg| = 2 \log{\sqrt{\frac{M-z}{\ubar{x}-z}} } + O(1) + O\left( \frac{\ubar{x}-z}{M-z}\right) \quad \text{as} \quad z \uparrow \ubar{x}.
\end{align*}
The integral $ \int_{z \vee \bar{x}}^M \frac{(T1) (s)}{ \sqrt{ (s-\bar{x})(s-z)}} \, ds$ can be handled analogously.

\end{proof}

\newpage

\subsubsection*{\hypertarget{algo L2 proj}{Minimization algorithm to find $Q^{V_n}$ projection of the naive estimator.}}
\scalebox{0.87}{
\begin{minipage}{\textwidth}
\begin{algorithm}[H] %or another one check
 \caption{Compute $V^{a^*}_n \equiv V^{\Pi}_n$}
    \DontPrintSemicolon
     \KwIn{vector observations $\mathbf{Z}$ and \textbf{Algorithm 2:} Compute $V^{a}_n$.}
    \KwOut{$V^{a^*}_n \equiv V^{\Pi}_n$}
     Using the preferred minimization procedure of the reader: \;
     $a^* = \argmin_{a \in \mathbb{R}^+} Q^{V_n}(V^a_n) = \argmin_{a \in \mathbb{R}^+} \int_{0}^{\infty} V^a_n(x)(V^a_n(x)-2V_n(x)) \, dx,$ \\
     where $V^a_n$ is computed using \textbf{Algorithm 2}.
\end{algorithm}

\begin{algorithm}[H]
\DontPrintSemicolon
\KwIn{$a>0$, vector observations $\mathbf{Z}$.}
\KwOut{$V^{a}_n$}

Compute $U_n$ and $U^*_n$, the least concave majorant of $U_n$.\;
Initialize an evaluation grid $\mathbf{x} = \left[x_1,\ldots,x_N\right]$.\;
\nlset{L-LCM} Compute the restricted Least Concave Majorant of $U_n$ from the Left, $\prescript{l}{}{U^*_n}$.\;
Initialize appropriately $\prescript{l}{}{U^*_n}$,  and $i=1$. Let $(U^*_n)^{\prime}$ be the r.h.s.\ derivative of $U^*_n$.\;
\vskip 0.1cm
\While{$ (U^*_n)^{\prime}(x_i)  >a$}{
$\prescript{l}{}{U^*_n}(x_i) \leftarrow U^*_n(x_i)$\;
$i \leftarrow i+1$\;
}
$i^* = \max \{ i \in \{1,\ldots,N \} \: : \:  (\prescript{l}{}{U}^*_n)^{\prime}(x_i) >a \}$\;
\For{$j \in \{i^*+1,\ldots, N\}$}{
$\prescript{l}{}{U^*_n}(x_i) \leftarrow a \cdot x_i$\;
}
\nlset{R-LCM} Compute the restricted Least Concave Majorant of $U_n$ from the Right, $\prescript{r}{}{U^*_n}$.\;
Initialize appropriately $\prescript{r}{}{U^*_n}$,  and $i=N$.\;
\vskip 0.1cm
\While{$(U^*_n)^{\prime}(x_i)  <a$}{
$\prescript{r}{}{U^*_n}(x_i) \leftarrow U^*_n(x_i)$\;
$i \leftarrow i-1$\;
}
$i_* = \min \{ i \in \{1,\ldots,N \} \: : \:   (\prescript{r}{}{U}^*_n)^{\prime}(x_i)  <a \}$\;
\For{$j \in \{1,\ldots, i_* -1\}$}{
$\prescript{r}{}{U^*_n}(x_i) \leftarrow a \cdot x_i$\;
}

$i_{\bar{x}} = \argmin_{i \: : \: x_i \geq \bar{x}}\{ |x_i - \bar{x}|\}$\;
$i_{\ubar{x}} = \argmin_{i \: : \: x_i \leq \ubar{x}}\{ |x_i - \ubar{x}|\}$\;

\If{$x_{i^*}< \ubar{x}$ or $x_{i^*}> \bar{x}$}{
$V^{a}_n = \begin{cases}
     (\prescript{l}{}{U}^*_n)^{\prime}(x_i) \quad \text{for} \: i \in \{1,\ldots,i_{\ubar{x}} -1\} \\
    a \quad \quad \quad \quad \quad \quad \quad \quad \quad \: \: \: \text{for} \: x_i \: : \: i \in \{i_{\ubar{x}},\ldots,i_{\bar{x}} \} \\
      (\prescript{r}{}{U}^*_n)^{\prime}(x_i) \quad \text{for} \: i \in \{i_{\bar{x}}+1,\ldots,N \}
\end{cases}$\;
}
\ElseIf{$x_{i^*} \in [\ubar{x},\bar{x}]$}{ 
$V^a_n = \begin{cases}
    ( U^*_n)^{\prime}(x_i)  \quad \text{for} \: i \in \{1,\ldots,i_{\ubar{x}}-1 \} \cup \{i_{\bar{x}}+1,\ldots,N \} \\
    a \quad \quad \quad \quad \quad \quad \quad \quad \quad \text{for} \: x_i \: : \: i \in \{i_{\ubar{x}},\ldots,i_{\bar{x}} \} 
\end{cases}$\;
}
\caption{Compute $V^{a}_n$.\label{IR}}
\end{algorithm}
\end{minipage}
}
\end{appendix}

\end{document}